\documentclass[a4paper,12pt]{article}

\usepackage{amsmath}
\usepackage{amsthm}
\usepackage{amssymb}
\usepackage{amscd}
\usepackage{hyperref}

\title{Estimating the eigenvalues on Quaternionic K\"ahler Manifolds}
\author{Yasushi Homma}
 \date{}

\theoremstyle{plain}
\newtheorem{theorem}{Theorem}[section]
\newtheorem{proposition}[theorem]{Proposition}
\newtheorem{corollary}[theorem]{Corollary}
\newtheorem{lemma}[theorem]{Lemma}
\theoremstyle{definition}

\numberwithin{equation}{section}

\theoremstyle{remark}
\newtheorem{remark}{Remark}[section]
\newtheorem{example}{Example}[section]

\newcommand{\End}{\mathrm{End}}
\newcommand{\id}{\mathrm{id}}
\newcommand{\sypc}{\mathfrak{sp}(n,\mathbb{C})}

\newcommand{\e}{\epsilon}
\newcommand{\sgn}{\mathrm{sign}}
\newcommand{\hotimes}{\hat{\otimes}}
\newcommand{\Sp}{\mathrm{Sp}}
\newcommand{\Cl}{\mathrm{Cl}}

\begin{document}

\maketitle

\begin{abstract}
We study geometric first order differential operators on quaternionic K\"ahler manifolds. Their principal symbols are related to the enveloping algebra and Casimir elements for $\Sp(1)\Sp(n)$. This observation leads to anti-symmetry of the principal symbols and Bochner-Weitzenb\"ock formulas for operators. As an application, we estimate the first eigenvalues of them.  
2000MSC: 58J60. 
17B35, 
53C26, 
\end{abstract}

\tableofcontents

\section{Introduction}\label{sec:1}
In differential geometry, Bochner-Weitzenb\"ock formulas play an important role to provide vanishing theorems and eigenvalue estimates for geometric differential operators. The strategy of giving such formulas is to find out algebraic structure among symbols of operators. As an example, we consider the Dirac operator on a spin manifold,
\[
D=\sum_i e_i\nabla_{e_i},
\]
where $\nabla$ is a covariant derivative and $\{e_i\}_i$ is a local orthonormal frame. The principal symbol of $D$ is the Clifford multiplication, which satisfies the Clifford relation $e_ie_j+e_je_i=-2\delta_{i,j}$. We rewrite this relation as
\begin{equation}
E_{ij}:=e_ie_j+\delta_{i,j}=-(e_je_i+\delta_{j,i})=-E_{ji}. \label{eqn:anti-symmetry}
\end{equation}
On the other hand, setting $\nabla^2_{X,Y}:=\nabla_X\nabla_Y-\nabla_{\nabla_XY}$, we know that the second order operator $\nabla^2_{e_i,e_j}$ satisfies a symmetric relation,
\begin{equation}
\nabla^2_{e_i,e_j}-R(e_i,e_j)/2=\nabla^2_{e_j,e_i}-R(e_j,e_i)/2,  \nonumber
\end{equation}
where $R$ is the curvature of $\nabla$. 
Combining this symmetry with the anti-symmetry \eqref{eqn:anti-symmetry}, we have
\[
\begin{split}
(D^2-\nabla^{\ast}\nabla-\kappa/4)=&\sum_{i,j} E_{ij}(\nabla^2_{e_i,e_j}-R(e_i,e_j)/2)\\
=&\sum_{i,j} -E_{ji}(\nabla^2_{e_j,e_i}-R(e_j,e_i)/2)
=-(D^2-\nabla^{\ast}\nabla-\kappa/4),
\end{split}
\]
and hence, $D^2=\nabla^{\ast}\nabla+\kappa/4$. Thus the essential point is to search anti-symmetry for the principal symbols of operators. 

As mentioned in \cite{H4}, the principal symbols of first order geometric operators called \textit{gradients} are controlled by the enveloping algebras of the Lie algebra $\mathfrak{so}(n)$ in a sense. From this observation, we have anti-symmetry such as \eqref{eqn:anti-symmetry} of the principal symbols.  Then we can give all Bochner-Weitzenb\"ock formulas for gradients on Riemannian or spin manifolds from a point of view of representation theory. In \cite{H3},  gradients on K\"ahler manifolds are also discussed. Working on the enveloping algebra associated to the holonomy group $\mathrm{U}(n)$, we can produce Bochner-Weitzenb\"ock formulas on K\"ahler manifolds. 

In this paper, we discuss gradients and their Bochner-Weitzenb\"ock formulas on  quaternionic K\"ahler manifolds with holonomy group in $\Sp(1)\Sp(n)$. The gradient is a first order differential operator defined to be an irreducible component of covariant derivative on an associated vector bundle. Most of differential operators in quaternionic K\"ahler geometry are realized as gradients or twisted gradients. 

In \cite{SW}, U. Semmelmann and G. Weingart presented an excellent method of giving  vanishing theorems and eigenvalue estimates for the Laplace operators and the Dirac operators on quaternionic K\"ahler manifolds. Their method is to consider twisted Dirac operators and compare the square of them with the Laplace operators. On the other hand, our method is more direct and universal. We give Bochner-Weitzenb\"ock formulas for gradients on any irreducible bundle concretely. Hence we can produce a lot of vanishing theorems and eigenvalue estimates. Thus our formulas are useful in various scenes of quaternionic geometry. As examples, we obtain lower bounds of the eigenvalues of the Laplace operators on differential forms. Our estimates not only cover the ones in \cite{SW} but are better than them.

This paper is organized as follows. In Section \ref{sec:2}, we give some results on the enveloping algebra of $\mathfrak{sp}(n)$. Proposition \ref{proposition:2-1} leads to anti-symmetry of symbols of gradients on quaternionic K\"ahler manifolds. In Section \ref{sec:3}, we calculate the eigenvalues of Casimir elements and state notation used in this paper. In Section \ref{sec:4}, we give anti-symmetry of the symbols as mentioned above. In Section \ref{sec:5}, we define gradients on quaternionic K\"ahler manifolds and show their conformal covariance, and in Section \ref{sec:6} we give the first main theorem, Bochner-Weitzenb\"ock formulas for gradients. After we discuss curvature endomorphisms on differential forms in Section \ref{sec:7}, we give the second main theorem, eigenvalue estimates of the Laplace operators, in the last section. Furthermore, we obtain some vanishing theorems, Proposition \ref{proposition:vanish} and Corollary \ref{corollary:vanish}.

Throughout this paper, we assume that the real dimension of a quaternionic K\"ahler manifold is greater than or equal to $8$. For the $4$-dimensional case, see \cite{H4}.

\section{Enveloping algebra and Casimir elements}\label{sec:2}

Let $E$ be a $2n$-dimensional complex vector space equipped with a complex symplectic structure $\sigma_E$, a quaternionic structure $J_E$, and a positive definite Hermitian inner product $\sigma_E(\cdot ,J_E(\cdot))$. We fix a symplectic unitary basis 
\[
\{\e_{\alpha}|\alpha=-n,-(n-1),\ldots,-1,1,\ldots,n\}
\]
such that
\begin{equation}
\sigma_E(\e_{\alpha},\e_{\beta})=\sgn(\alpha)\delta_{\alpha,-\beta},\quad J_E(\e_{\alpha})=\sgn(\alpha)\e_{-\alpha}, \nonumber
\end{equation}
where $\sgn(\alpha)=\pm 1$ is the sign of $\alpha$. 

The complex symplectic group $\Sp(n,\mathbb{C})$ on $E$ is the group of automorphisms preserving $\sigma_E$ and the symplectic group $\Sp(n)$ is the real subgroup of $\Sp(n,\mathbb{C})$ compatible with $J_E$. The second symmetric tensor product space $S^2(E)$ is isomorphic to the Lie algebra $\mathfrak{sp}(n,\mathbb{C})$ of $\Sp(n,\mathbb{C})$ by associating $\e\odot \e'$ in $S^2(E)$ with an endomorphism
\begin{equation}
(\e\odot \e')(u):=\sigma_E(\e,u)\e'+\sigma_E(\e',u)\e \quad \textrm{for each $u$ in $E$}.\nonumber
\end{equation}
In particular, we can choose $\{\e_{\alpha}\odot \e_{\beta}|\alpha+\beta \ge 0\}$ as a basis of $\sypc$. But, to construct Casimir elements, it is better to employ $x_{\alpha\beta}:=-\sgn(\beta)\e_{\alpha}\odot \e_{-\beta}$ instead of $\e_{\alpha}\odot \e_{\beta}$ as a basis. This $x_{\alpha\beta}$ acts on $E$ by
\[
x_{\alpha\beta}(\e_{\nu})=\delta_{\beta,\nu}\e_{\alpha}-\sgn(\alpha\beta)\delta_{-\alpha,\nu}\e_{-\beta},
\]
 and satisfies 
\begin{gather}
x_{\alpha\beta}=-\sgn(\alpha\beta)x_{-\beta-\alpha},\nonumber\\ 
[x_{\alpha\beta},x_{\mu\nu}]=\delta_{\beta,\mu}x_{\alpha\nu}-\delta_{\alpha,\nu}x_{\mu\beta}+\sgn(\alpha\beta)(\delta_{-\beta,\nu}x_{\mu-\alpha}-\delta_{-\alpha,\mu}x_{-\beta\nu})\nonumber
\end{gather}
for $\alpha,\beta,\mu,\nu=\pm 1,\dots,\pm n$. 

Let $U(\sypc)$ be the universal enveloping algebra of $\sypc$. The center $\mathfrak{Z}$ of $U(\sypc)$ is characterized as the invariant sub-algebra of $U(\sypc)$ under the adjoint action of $\Sp(n,\mathbb{C})$, whose elements are called \textit{Casimir elements}. It is well-known how to construct Casimir elements generating $\mathfrak{Z}$ \cite{O}, \cite{Z}. For each nonnegative integer $q$, we define an element $x_{\alpha\beta}^q$ in $U(\sypc)$ by
\begin{equation}
x_{\alpha\beta}^q:=
\begin{cases}
\sum_{\alpha_1,\alpha_2,\dots,\alpha_{q-1}=\pm 1,\dots,\pm n}x_{\alpha\alpha_1}x_{\alpha_1\alpha_2}\cdots x_{\alpha_{q-1}\beta} & q\ge 1\\
\delta_{\alpha,\beta}  & q=0.
\end{cases}\nonumber
\end{equation}
Then the trace $c_q:=\sum_{\alpha} x_{\alpha\alpha}^q$ is a Casimir element and the center $\mathfrak{Z}$ is generated by $c_2,c_4,\ldots,c_{2n}$. We will need their translated elements in later sections. The translated elements  are defined by 
\begin{gather}
\hat{x}_{\alpha\beta}:=x_{\alpha\beta}-(n+1/2)\delta_{\alpha,\beta},\nonumber\\
\hat{x}_{\alpha\beta}^q:=
\begin{cases}
\sum_{\alpha_1,\alpha_2,\dots,\alpha_{q-1}=\pm 1,\dots,\pm n}\hat{x}_{\alpha\alpha_1}\hat{x}_{\alpha_1\alpha_2}\cdots \hat{x}_{\alpha_{q-1}\beta} & q\ge 1\\
\delta_{\alpha,\beta}  & q=0,
\end{cases}\nonumber\\
\hat{c}_q:=\sum_{\alpha} \hat{x}_{\alpha\alpha}^q. \nonumber
\end{gather}
In \cite{H2}, the author showed the following.
\begin{proposition}\label{proposition:2-1}
Translated elements $\{\hat{x}_{\alpha\beta}^q\}_{q,\alpha,\beta}$ satisfy
\begin{equation}
\hat{x}_{\alpha\beta}^q=\sgn(\alpha\beta)\left\{(-1)^q\hat{x}^q_{-\beta-\alpha}-\frac{1-(-1)^q}{2}\hat{x}_{-\beta-\alpha}^{q-1}-\sum_{p=0}^{q-1}(-1)^p\hat{c}_{q-1-p}\hat{x}_{-\beta-\alpha}^p\right\},\label{eqn:ubw}
\end{equation}
and
\begin{equation}
\hat{x}_{\alpha\beta}^q=\sum_{p=0}^q\binom{q}{p}\left(-n-\frac{1}{2}\right)^{q-p}x_{\alpha\beta}^p.\label{eqn:hat}
\end{equation}
The translated Casimir elements $\{\hat{c}_q\}_q$ satisfy 
\begin{equation}
2\hat{c}_{2q+1}=-\hat{c}_{2q}-\sum_{p=0}^{2q}(-1)^p\hat{c}_{2q-p}\hat{c}_p.
\label{eqn:casimir-1}
\end{equation}
\end{proposition}
\section{Representation of $\Sp(n)$ and eigenvalues of Casimir elements}\label{sec:3}
We set $\mathfrak{h}:=\mathrm{span}_{\mathbb{R}}\{\sqrt{-1}x_{ii}|i=1,\dots,n\}$ as a maximal abelian subalgebra of $\mathfrak{sp}(n)$. We consider a finite-dimensional irreducible unitary $\Sp(n)$-module $V$ and decompose it into simultaneous eigenspaces with respect to $\mathfrak{h}$, $V=\bigoplus_{\lambda}V(\lambda)$. Each eigenvalue $\lambda=(\lambda^1,\dots,\lambda^n)$ is called weight, and a weight vector $\phi_{\lambda}$ in $V(\lambda)$ satisfies $x_{ii}\phi_{\lambda}=\lambda^i \phi_{\lambda}$ for $i=1,\dots,n$. Ordering these weights lexicographically, we have the highest weight $\rho$ which satisfies the dominant integral condition,
\begin{equation}
\rho=(\rho^1,\dots,\rho^n)\in \mathbb{Z}^n \textrm{ and } \rho^1\ge \rho^2\ge\cdots \ge \rho^n\ge 0. \nonumber
\end{equation}
Conversely, for a dominant integral weight $\rho$, we can construct an irreducible unitary $\Sp(n)$-module with highest weight $\rho$.  Therefore we denote by $(\pi_{\rho},V_{\rho})$ an irreducible unitary representation of $\Sp(n)$ and its infinitesimal one with highest weight $\rho$. When writing a weight, we denote by $k_a$ a string of $k$ with length $a$ and sometimes omit a string of $0$. For example, $(1_a)=(1_a,0_{n-a})$ is a weight such that the first $a$ components are $1$ and the others are $0$. 
\begin{example}
The quaternionic vector space $E$ is an irreducible $\Sp(n)$-module with highest weight $(1_1)$. The second symmetric tensor product space $S^2(E)\simeq \mathfrak{sp}(n,\mathbb{C})$ has the highest weight $(2_1)$. The exterior tensor product space $\Lambda^a(E)$ for $2\le a\le n$ is not irreducible. Actually, by using the complex symplectic form $\sigma_E$, we can decompose $\Lambda^a(E)$ into irreducible components, 
\begin{equation}
\Lambda^a(E)=\bigoplus_{p=0}^{[a/2]}\sigma_E^{p}\Lambda^{a-2p}_0(E). \nonumber
\end{equation}
Here, $\Lambda^{a}_0(E)$ is the so-called primitive space of $\Lambda^a(E)$, which is an irreducible $\Sp(n)$-module with highest weight $(1_a)$.
\end{example}
\begin{example}
In quaternionic K\"ahler geometry, we often discuss an $\Sp(n)$-module with highest weight $(2_b,1_{a-b})$ for $0\le b \le a\le n$. The representation space is realized as the top irreducible summand of $\Lambda^a_0(E)\otimes \Lambda^b_0(E)$. We denote it by $\Lambda^{a,b}_0(E)$. 
\end{example}
\begin{remark}
By  quaternionic structure $J_E$ on $E$, we can put a quaternionic or real structure on each irreducible $\Sp(n)$-module \cite{H2}. In fact, there is a quaternionic (resp. real) structure on $V_{\rho}$ in the case that $\sum_i \rho^i$ is odd (resp. even). 
\end{remark}

We shall calculate eigenvalues of Casimir element $c_q$ on irreducible $\Sp(n)$-modules. We consider a tensor product space $V_{\rho}\otimes E$. The highest weights of irreducible summands in $V_{\rho}\otimes E$ are
\begin{equation}
\{\rho+\mu_{\nu}\:|\: \textrm{$\rho+\mu_{\nu}$ is dominant integral}, \nu=\pm 1,\dots,\pm n \}, \nonumber
\end{equation}
where
\begin{equation}
\mu_{\nu}=
\begin{cases}
\mu_{\nu}:=(\underbrace{0,\dots,0}_{\nu-1},1,\underbrace{0,\dots,0}_{n-\nu}) & \textrm{for $1\le \nu \le n$},\\
\mu_{\nu}:=-\mu_{-\nu}& \textrm{for $-n\le \nu \le -1$}.
\end{cases}\nonumber
\end{equation}
Setting $V_{\rho+\mu_{\nu}}:=\{0\}$ for $\rho+\mu_{\nu}$ without dominant integral condition, we can describe the irreducible decomposition of $V_{\rho}\otimes E$ as 
\begin{equation}
V_{\rho}\otimes E=\bigoplus_{\nu=\pm 1,\dots,\pm n}V_{\rho+\mu_{\nu}}=\bigoplus_{i=1,\dots,n}(V_{\rho+\mu_{i}}\oplus V_{\rho-\mu_i}). \nonumber
\end{equation}
Each component $V_{\rho+\mu_{\nu}}$ is equipped with  a Hermitian inner product of the restriction of the one on $V_{\rho}\otimes E$. 

Let $\Pi_{\nu}$ be the orthogonal projection from $V_{\rho}\otimes E$ onto $V_{\rho+\mu_{\nu}}$. We define a linear mapping $p_{\nu}(\e):V_{\rho}\to V_{\rho+\mu_{\nu}}$ for each $\e$ in $E$ by
\begin{equation}
V_{\rho}\ni \phi\mapsto p_{\nu}(\e)\phi:=\Pi_{\nu}(\phi\otimes \e)\in V_{\rho+\mu_{\nu}},\label{proj}
\end{equation}
and denote by $p_{\nu}(\e)^{\ast}$ the adjoint map of $p_{\nu}(\e)$ with respect to Hermitian inner products of $V_{\rho}$ and $V_{\rho+\mu_{\nu}}$. To connect with these linear mappings and the enveloping algebra, we assign a constant $w_{\nu}$ called \textit{the conformal weight} to $\rho+\mu_{\nu}$, 
\begin{equation}
\begin{cases}
w_{i}:=-(\rho^i-i+1)   & \textrm{to $\rho+\mu_i$ for $i=1,\dots,n$}, \\
w_{-i}:=\rho^i-i+2n+1  & \textrm{to $\rho+\mu_{-i}$ for $i=1,\dots,n$},
\end{cases}\nonumber
\end{equation}
and define the translated conformal weight $\hat{w}_{\nu}$ by 
\begin{equation}
\hat{w}_{\nu}:=w_{\nu}-(n+1/2). \nonumber
\end{equation}
Then we have 
\begin{proposition}[\cite{H2},\cite{O}]\label{proposition:3-1}
\begin{gather}
\sum_{\nu=\pm 1,\dots, \pm n} w_{\nu}^q p_{\nu} (\e_{\alpha})^{\ast} 
  p_{\nu}(\e_{\beta})=\sgn(\alpha\beta)\pi_{\rho}(x_{-\alpha-\beta}^q), \label{eqn:cl}\\
\sum_{\nu=\pm 1,\dots, \pm n}\hat{w}_{\nu}^qp_{\nu}(\e_{\alpha})^{\ast}
 p_{\nu}(\e_{\beta})=\sgn(\alpha\beta)\pi_{\rho}(\hat{x}_{-\alpha-\beta}^q), 
  \label{eqn:cl-2} \\
 \sum_{\alpha=\pm 1,\dots,\pm n} p_{\nu}(\e_{\alpha})^{\ast}p_{\nu}(\e_{\alpha})=\frac{ \dim V_{\rho+\mu_{\nu}} }{\dim V_{\rho}}, \nonumber \\
\pi_{\rho}(c_q)=\sum_{\nu=\pm 1,\dots, \pm n} w_{\nu}^q \frac{\dim V_{\rho+\mu_{\nu}}}{\dim V_{\rho}},\quad \pi_{\rho}(\hat{c}_q)=\sum_{\nu=\pm 1,\dots, \pm n} \hat{w}_{\nu}^q \frac{\dim V_{\rho+\mu_{\nu}}}{\dim V_{\rho}}. \label{eqn:casimir-2}
\end{gather}
\end{proposition}
Without this proposition, we can compute the eigenvalues of $c_q$ and $\hat{c}_q$ for $q=0,1,2$ \cite{O}, \cite{Z}. On an irreducible $\Sp(n)$-module $V_{\rho}$, 
\begin{equation}
\begin{split}
\pi_{\rho}(c_0)=2n,\quad \pi_{\rho}(c_1)=0, \quad \pi_{\rho}(c_2)=2\sum_{i=1}^n \rho^i\left(\rho^i+2(n-i+1)\right), \\
\pi_{\rho}(\hat{c}_0)=2n, \quad \pi_{\rho}(\hat{c}_1)=-2n^2-n,\quad \pi_{\rho}(\hat{c}_2)=\pi_{\rho}(c_2)+2n(n+1/2)^2.
\end{split}\label{eqn:c-1}
\end{equation}
From \eqref{eqn:casimir-1}, we also have the eigenvalue of $c_3$ and $\hat{c}_3$, 
\begin{equation}
\pi_{\rho}(c_3)=(n+1)\pi_{\rho}(c_2), \quad \pi_{\rho}(\hat{c}_3)=-(2n+1/2)\pi_{\rho}(c_2)-2n(n+1/2)^3. \label{eqn:c-2}
\end{equation}
To calculate the eigenvalues of higher Casimir elements, we need Proposition \ref{proposition:3-1}. In fact, there are formulas for the eigenvalues of higher Casimir elements \cite{O}, \cite{Z}. But those formulas are complicated to compute explicitly. On the other hand, D. Calderbank, P. Gauduchon and M. Herzlich gave a nice formula to calculate the eigenvalues of Casimir elements for the special orthogonal group $\mathrm{SO}(n)$ \cite{CGH}. From a similar discussion, we have a formula of $\pi_{\rho}(c_{q})$ calculated easily. 
\begin{proposition}\label{proposition:3-2}
We denote by $\mathcal{N}$ the number of irreducible summands in $V_{\rho}\otimes E$. Then the relative dimension $\dim V_{\rho+\mu_{\nu}}/\dim V_{\rho}$ is given by
\begin{equation}
\frac{\dim V_{\rho+\mu_{\nu}}}{\dim V_{\rho}}=-2(\hat{w}_{\nu}-(-1)^{\mathcal{N}}) \prod_{\scriptstyle \nu' \neq \nu, \atop \scriptstyle    \textrm{ $\rho+\mu_{\nu'}$ is dominant} }  \frac{\hat{w}_{\nu}+\hat{w}_{\nu'}}{\hat{w}_{\nu}-\hat{w}_{\nu'}}.
\nonumber
\end{equation}
Of course, the relative dimension $\dim V_{\rho+\mu_{\nu}}/\dim V_{\rho}$ is zero for $\rho+\mu_{\nu}$ without dominant integral condition. By the above equation and \eqref{eqn:casimir-2}, we can easily calculate the eigenvalues of $c_q$.
\end{proposition}
\begin{remark}
The number $\mathcal{N}$ is odd if and only if the $n$th component $\rho^n$ of the highest weight $\rho$ is zero.
\end{remark}
\begin{example}
We shall calculate eigenvalues of $c_2$ and $c_4$ on $V_{(2_b,1_{a-b})}=\Lambda^{a,b}_0(E)$. The $\Sp(n)$-module $V_{(2_b,1_{a-b})}\otimes E$ splits as 
\[
\begin{split}
V_{\rho}\otimes E=&V_{\rho+\mu_1}\oplus V_{\rho+\mu_{b+1}}\oplus V_{\rho+\mu_{a+1}}\oplus V_{\rho+\mu_{-b}}\oplus V_{\rho+\mu_{-a}}\\
 =&V_{(3,2_{b-1},1_{a-b})}\oplus V_{(2_{b+1},1_{a-b-1})}\oplus V_{(2_{b},1_{a-b+1})}\oplus V_{(2_{b-1},1_{a-b+1})}\oplus V_{(2_{b},1_{a-b-1})}.
\end{split}
\]
The next table of the relative dimension $\dim V_{\rho+\mu_{\nu}}/\dim V_{\rho}$ follows from Proposition \ref{proposition:3-2}.   
\begin{table}[h]
\begin{tabular}{|c|c|c|}
  \hline  
$\rho+\mu_{\nu}$ &  $w_{\nu}$    & relative dimension \\
\hline
 $\rho+\mu_1$ &   $-2$  & $\displaystyle{\frac{2b(a+1)(2n-a+3)(2n-b+4)(n+2)  }{(a+2)(b+1)(2n-a+4)(2n-b+5)}}$\\
\hline
$\rho+\mu_{b+1}$ &   $b-1$  &$\displaystyle{\frac{(a-b)(2n-b+4)(2n-a-b+2)(n-b+1)}{(b+1)(a-b+1)(2n-a-b+3)(n-b+2)}}$\\
\hline
$\rho+\mu_{a+1}$ &   $a$    & $\displaystyle{\frac{(a-b+2)(2n-a+3)(2n-a-b+2)(n-a)}{(a+2)(a-b+1)(2n-a-b+3)(n-a+1)}}$\\
\hline
$\rho+\mu_{-b}$ &   $2n-b+3$  & $\displaystyle{\frac{b(a-b+2)(2n-a-b+4)(n-b+3)}{(a-b+1)(2n-b+5)(2n-a-b+3)(n-b+2)}}$\\
\hline
$\rho+\mu_{-a}$ &   $2n-a+2$  & $\displaystyle{\frac{(a+1)(a-b)(2n-a-b+4)(n-a+2)}{(a-b+1)(2n-a+4)(2n-a-b+3)(n-a+1)}}$\\
\hline
\end{tabular}
\caption{}\label{table-1}
\end{table}%

Then we have 
\begin{equation}
\begin{split}
\pi_{(2_b,1_{a-b})}(c_2)=&  2a(2n-a+2)+2b(2n-b+4),\\
\pi_{(2_b,1_{a-b})}(c_4)=&2a(2n-a+2)(2n+3)(n+1) -2a^2(2n-a+2)^2\\
   &+2b(2n-b+4)(2n+3)(n+3)-2b^2(2n-b+4)^2.
\end{split}\nonumber
\end{equation}
\end{example}

In the next section, we will discuss symbols of gradients on a quaternionic K\"ahler manifold whose holonomy group is in $\Sp(1)\Sp(n)$. Then we shall state some facts for the $\Sp(1)$-case. We consider a $2$-dimensional complex vector space $H$ with quaternionic structure $J_H$ and symplectic structure $\sigma_H$. The group of automorphisms on $H$ preserving $J_H$ and $\sigma_H$ is the Lie group $\Sp(1)$. In other words, $(H,J_H,\sigma_H)$ is the natural $\Sp(1)$-module. We set $\{h_A\}_{A=\pm 1}$ as a symplectic unitary basis of $H$ and 
\[
\{y_{AB}:=-\sgn(B)h_{A}\odot h_{-B}\;|\; A,B=\pm 1,A+B\ge 0\}
\] as a basis of $\mathfrak{sp}(1,\mathbb{C})\simeq S^2(H)$. We will use only the following elements and relations in the enveloping algebra $U(\mathfrak{sp}(1,\mathbb{C}))$,
\begin{equation}
y_{AB}^0:=\delta_{A,B},\quad y_{AB}^1:=y_{AB}, \quad C_2:=\sum_{A,B}y_{AB}y_{BA} \nonumber
\end{equation}
and
\begin{equation}
y_{AB}=-\sgn(AB)y_{-B-A}\quad  \textrm{for $A,B=\pm 1$}.\label{eqn:ubw-2}
\end{equation}
Since all the irreducible $\Sp(1)$-modules are parametrized by non-negative integer $k$, we denote by $(\pi_k,V_k)$ an irreducible $\Sp(1)$-module with highest weight $k$. Note that $V_k$ is isomorphic to the $k$th symmetric tensor product space $S^k(H)$ of $H$, and $\pi_k(C_2)$ is $2k(k+2)$. 

We consider $V_k\otimes H$ and decompose it, $V_k\otimes H=V_{k+1}\oplus V_{k-1}$. In the same way as \eqref{proj}, we define a linear map $p_{N}(\cdot)$ from $V_k$ to $V_{k+N}$ for $N=\pm 1$ and denote by $W_N$ the conformal weight associated to $k+N$. Here, $W_{1}=-k$ and $W_{-1}=k+2$. The equation \eqref{eqn:cl} for the $\Sp(1)$-case is 
\begin{equation}
\sum_{N\pm 1} W_N^q p_N(h_A)^{\ast}p_N(h_B)=\sgn(AB)\pi_k(y^q_{-A-B}). \label{eqn:cl-3}
\end{equation}

We summarize notation used in this paper.
\begin{table}[h]
\begin{tabular}{|c|c|c|}
\hline
   &    $G=\Sp(n)$  & $G=\Sp(1)$ \\ 
\hline 
suffices & $\alpha,\beta, \nu =\pm 1,\dots,\pm n$ & $A,B,N=\pm 1$\\   
\hline
the natural $G$-module & $(E,J_E,\sigma_E) $ & $(H,J_H,\sigma_H) $\\
\hline
a (symplectic) unitary basis & $\{ \e_{\alpha}\}_{\alpha}$ & $\{ h_{A} \}_A$ \\
\hline 
 a basis of the Lie algebra & $\{x_{\alpha\beta}|\alpha+\beta\ge 0\}$ & $\{y_{AB}|A+B\ge 0\}$ \\
\hline
Casimir elements  &$ \{c_q\}_{q\ge 0}$,  $\{\hat{c}_q\}_{q\ge 0}$  & $C_2$ \\
\hline
irreducible $G$-module & $(\pi_{\rho},V_{\rho})$ & $(\pi_{k},V_{k})$\\
\hline
the dominant integral & $\rho=(\rho^1,\dots,\rho^n)\in \mathbb{Z}^n$, & $k\in \mathbb{Z}$, \\
         condition  &  $\rho^1\ge \cdots \ge \rho^n\ge 0$  & $k\ge 0$   \\
\hline 
linear mapping & $p_{\nu}(\cdot):V_{\rho}\to V_{\rho+\mu_{\nu}}$ & $p_N(\cdot):V_k\to V_{k+N}$ \\
\hline
conformal weights & $w_{\nu}$,  $\hat{w}_{\nu}$    & $W_N$ \\
\hline

\end{tabular}
\caption[]{}
\label{table-2}
\end{table}%

\section{Symbols of gradients and their relations}\label{sec:4}
Let $(H,J_H,\sigma_H)$ (resp. $(E,J_E,\sigma_E)$) be the natural $\Sp(1)$-module (resp. $\Sp(n)$-module). The outer tensor product space $H\hotimes E$ has a real structure $J_H\hotimes J_E$ and a Hermitian structure $\sigma_H(\cdot,J_H(\cdot))\hotimes \sigma_E(\cdot,J_E(\cdot))$. The real part of this vector space is a real $4n$-dimensional vector space with positive inner product. This $H\hotimes E$ is a model of tangent space of a quaternionic K\"ahler manifold. Taking unitary bases $\{h_A\}_{A}$ of $H$ and $\{\e_{\alpha}\}_{\alpha}$ of $E$, we have the one of $H\hotimes E$,
\begin{equation}
\{v_{A,\alpha}:=h_A\hotimes \e_{\alpha} \;|\; A=\pm 1,\alpha=\pm 1,\dots,\pm n\}. \nonumber 
\end{equation}

We consider the Lie groups $\Sp(1)\times \Sp(n)$ and $\Sp(1)\Sp(n):=(\Sp(1)\times \Sp(n))/\{\pm I\}$. Each unitary irreducible representation of $\Sp(1)\times \Sp(n)$ is given by  
\[
(\pi_{k,\rho},V_{k,\rho}):=(\pi_k\hotimes \pi_{\rho},V_k\hotimes V_{\rho}).
\]
When $k+\sum_{i} \rho^i$ is odd, $(\pi_{k,\rho},V_{k,\rho})$ does not factor through a representation of $\Sp(1)\Sp(n)$. Furthermore, from quaternionic or real structures on $V_k$ and $V_{\rho}$, we can set a quaternionic structure on $V_{k,\rho}$. When $k+\sum_i \rho^i$ is even, $V_{k,\rho}$ is an irreducible $\Sp(1)\Sp(n)$-module with real structure. For example, $H\hotimes E=V_1\hotimes V_{(1_1)}$ is an irreducible $\Sp(1)\Sp(n)$-module with real structure $J_H\hotimes J_E$. 

Now, we consider a representation space $V_{k,\rho}\otimes(H\hotimes E)$ of $\Sp(1)\times \Sp(n)$ or $\Sp(1)\Sp(n)$ and decompose it,
\begin{equation}
V_{k,\rho}\otimes (H\hotimes E)=\bigoplus_{\scriptstyle N=\pm 1, \atop \scriptstyle \nu=\pm 1,\dots,\pm n}V_{k+N,\rho+\mu_{\nu}}. \nonumber
\end{equation}
Considering the orthogonal projection from $V_{k,\rho}\otimes(H\hotimes E)$ onto $V_{k+N,\rho+\mu_{\nu}}$, we define a linear mapping $p_{N,\nu}(\cdot)$ from $V_{k,\rho}$ to $V_{k+N,\rho+\mu_{\nu}}$ by
\begin{equation}
p_{N,\nu}(h\hotimes \e)(\phi\hotimes \psi):=p_N(h)\phi\hotimes p_{\nu}(\e)\psi
\nonumber
\end{equation}
for $h\hotimes \e$ in $H\hotimes E$ and $\phi\hotimes \psi$ in $V_{k,\rho}=V_k\hotimes V_{\rho}$. The adjoint map of $p_{N,\nu}(h\hotimes \e)$ is defined to be $p_{N}(h)^{\ast}\hotimes p_{\nu}(\e)^{\ast}$. These maps are just the principal symbols of first order differential operators called quaternionic K\"ahlerian gradients given in the next section. As mentioned in Section \ref{sec:1}, to obtain Bochner-Weitzenb\"ock formulas for the operators, we need anti-symmetric relations among the principal symbols. Such relations follow from the next proposition. 
\begin{proposition}\label{proposition:4-1}
The linear maps $\{p_{N,\nu}(v_{A,\alpha})\}_{N,\nu,A,\alpha}$ on $V_{k,\rho}$ satisfy the followings. 
\begin{enumerate}
\item
\begin{equation}
\begin{split}
 &\sum_{N, \nu}\hat{w}_{\nu}^qp_{N,\nu}(v_{A,\alpha})^{\ast}p_{N,\nu}(v_{B,\beta}) \\
=& \sgn (AB \alpha\beta)\sum_{N, \nu}\left\{(-1)^q\hat{w}_{\nu}^q-\frac{1-(-1)^q}{2}\hat{w}_{\nu}^{q-1} \right.\\
&\left.-\sum_{p=0}^{q-1}(-1)^p \pi_{\rho}(\hat{c}_{q-1-p})\hat{w}_{\nu}^p\right\}p_{N,\nu}(v_{-B,-\beta})^{\ast}p_{N,\nu}(v_{-A,-\alpha}).
\end{split} \label{relation-1}
\end{equation}
\item When $k$ is not zero, 
\begin{equation}
\begin{split}
&\sum_{N,\nu}W_N\hat{w}_{\nu}^qp_{N,\nu}(v_{A,\alpha})^{\ast}p_{N,\nu}(v_{B,\beta})\\
=& -\sgn (AB \alpha\beta)\sum_{N,\nu}W_N\left\{(-1)^q\hat{w}_{\nu}^q-\frac{1-(-1)^q}{2}\hat{w}_{\nu}^{q-1}\right.
\\
&\left.-\sum_{p=0}^{q-1}(-1)^p \pi_{\rho}(\hat{c}_{q-1-p})\hat{w}_{\nu}^p\right\}
p_{N,\nu}(v_{-B,-\beta})^{\ast}p_{N,\nu}(v_{-A,-\alpha}).
\end{split}\label{relation-2}
\end{equation}
\end{enumerate}
\end{proposition}
\begin{proof}
It follows from \eqref{eqn:cl-2} and \eqref{eqn:cl-3} that 
\begin{gather}
\sum_{N,\nu}\hat{w}_{\nu}^qp_{N,\nu}(v_{A,\alpha})^{\ast}p_{N,\nu}(v_{B,\beta})=\sgn(AB\alpha\beta)\delta_{A,B} \id \hotimes \pi_{\rho}(\hat{x}^q_{-\alpha-\beta}), 
\label{eqn:4-4}\\
\sum_{N,\nu}W_N\hat{w}_{\nu}^qp_{N,\nu}(v_{A,\alpha})^{\ast}p_{N,\nu}(v_{B,\beta})=\sgn(AB\alpha\beta)\pi_k(y_{-A-B})\hotimes \pi_{\rho}(\hat{x}^q_{-\alpha-\beta}). \nonumber
\end{gather}
Substituting \eqref{eqn:ubw} and \eqref{eqn:ubw-2} for the above equations, we can prove the proposition. 
\end{proof}

\section{Gradients on quaternionic K\"ahler manifolds }\label{sec:5}
Let $(M,g)$ be a real  $4n$-dimensional quaternionic K\"ahler manifold. The frame bundle of $M$ reduces to a principal bundle $\mathbf{P}$ with structure group $\Sp(1)\Sp(n)$. Take the $\Sp(1)\Sp(n)$-module $H\hotimes E$, and we have an associated vector bundle $\mathbf{H}\hotimes \mathbf{E}:=\mathbf{P}\times_{\Sp(1)\Sp(n)}(H\hotimes E)$ with real structure and Hermitian metric. The real part of $\mathbf{H}\hotimes \mathbf{E}$ is isometric to the tangent bundle $T(M)$.

For an irreducible $\Sp(1)\Sp(n)$-module $V_{k,\rho}$, we have an associated vector bundle $\mathbf{S}_{k,\rho}:=\mathbf{P}\times_{\Sp(1)\Sp(n)}V_{k,\rho}$. Since the Levi-Civita connection reduces a connection on $\mathbf{P}$, we have a covariant derivative $\nabla$ on $\mathbf{S}_{k,\rho}$,
\begin{equation}
\nabla:\Gamma(M,\mathbf{S}_{k,\rho})\to \Gamma (M,\mathbf{S}_{k,\rho}\otimes (T^{\ast}(M)\otimes \mathbb{C}))\simeq \Gamma (M,\mathbf{S}_{k,\rho}\otimes(\mathbf{H}\hotimes \mathbf{E})). \nonumber
\end{equation}
Here we identified $T^{\ast}(M)\otimes \mathbb{C}$ with $T(M)\otimes \mathbb{C}\simeq \mathbf{H}\hotimes \mathbf{E}$ by complex inner product $\sigma_H\hotimes \sigma_E$. Decomposing $\nabla$ with respect to $\Sp(1)\Sp(n)$, we have first order differential operators in the following way. Let $\{h_A\}_A$ and $\{\e_{\alpha}\}_{\alpha}$ be local unitary frames  of $\mathbf{H}$ and $\mathbf{E}$, respectively. For a smooth section $\phi$ of $\mathbf{S}_{k,\rho}$, the derivative $\nabla\phi$ is locally expressed by 
\begin{equation}
\nabla\phi=\sum_{A,\alpha} \nabla_{v_{A,\alpha}}\phi \otimes v_{A,\alpha}^{\ast}=\sum_{A,\alpha} \sgn(A\alpha)\nabla_{v_{A,\alpha}}\phi \otimes v_{-A,-\alpha},
\nonumber
\end{equation}
where we set $v_{A,\alpha}:=h_A\hotimes \e_{\alpha}$. We project $\nabla \phi$ from $\mathbf{S}_{k,\rho}\otimes(\mathbf{H}\hotimes \mathbf{E})$ onto an irreducible bundle $\mathbf{S}_{k+N,\rho+\mu_{\nu}}$ fiberwise. Then we define a first order differential operator $D_{N,\nu}:\Gamma(M,\mathbf{S}_{k,\rho})\to \Gamma(M,\mathbf{S}_{k+N,\rho+\mu_{\nu}})$ by 
\begin{equation}
D_{N,\nu}:=\sum_{A,\alpha} \sgn(A\alpha)p_{N,\nu}(v_{-A,-\alpha})\nabla_{v_{A,\alpha}}. \label{eqn:gradient}
\end{equation}
It is easy to show that the formal adjoint operator $(D_{N\nu})^{\ast}$ of $D_{N\nu}$ is given by
\begin{equation}
(D_{N,\nu})^{\ast}=-\sum p_{N,\nu}(v_{A,\alpha})^{\ast}\nabla_{v_{A,\alpha}} :\Gamma(M,\mathbf{S}_{k+N,\rho+\mu_{\nu}})\to \Gamma(M,\mathbf{S}_{k,\rho}).
\nonumber
\end{equation}
If $k+N$ or $\rho+\mu_{\nu}$ does not satisfy dominant integral condition, we set $\mathbf{S}_{k+N,\rho+\mu_{\nu}}:=M\times\{0\}$ and $D_{N,\nu}:=0$ virtually.  We call these operators $\{D_{N,\nu},(D_{N,\nu})^{\ast}\}_{N,\nu}$ \textit{gradients on a quaternionic K\"ahler manifold} or \textit{quaternionic K\"ahlerian gradients}. Most of first order differential operators in quaternionic K\"ahler geometry are realized as gradients. 
\begin{remark}
The obstruction for lifting $\mathbf{P}$ to a principal $\Sp(1)\times \Sp(n)$ bundle $\tilde{\mathbf{P}}$ is the second Stiefel-Whitney class of the real part of $S^2(\mathbf{H})$ (cf. \cite{S1}). When the obstruction is zero, we can consider a vector bundle $\mathbf{S}_{k,\rho}$ associated to $\tilde{\mathbf{P}}$. In this case, we have gradients even if $k+\sum \rho^i$ is odd. Since our results come from local calculation, we can do well in the case that $k+\sum \rho^i$ is odd. 
\end{remark}

\begin{example}
Let $M$ be a spin manifold and $\mathbf{S}(M)$ be the spinor bundle. The Dirac operator $D$ is realized as an irreducible  component of $\nabla$ on $\mathbf{S}(M)$,
\begin{equation}
D:\Gamma(M,\mathbf{S}(M))\xrightarrow{\nabla}\Gamma(M,\mathbf{S}(M)\otimes (T^{\ast}(M)\otimes\mathbb{C}))\xrightarrow{\Pi} \Gamma(M,\mathbf{S}(M)). \nonumber
\end{equation}
If $M$ has a quaternionic K\"ahler structure, then the spinor bundle $\mathbf{S}(M)$ is decomposed with respect to $\Sp(1)\Sp(n)$ or $\Sp(1)\times \Sp(n)$,
\begin{equation}
\mathbf{S}(M)=\bigoplus_{k=0}^n \mathbf{S}_{k,(1_{n-k})}=\bigoplus_{k=0}^n S^k(\mathbf{H})\hotimes \Lambda^{n-k}_0(\mathbf{E}). \nonumber
\end{equation}
We divide the Dirac operator $D$ along this decomposition. Then each piece of $D$ is a quaternionic K\"ahlerian gradient (cf. \cite{KSW}). 
\end{example}

In the rest of this section, we show a conformal covariance of gradient. An almost quaternionic Hermitian manifold $(M,g)$ is a $4n$-dimensional Riemannian manifold whose frame bundle reduces to an $\Sp(1)\Sp(n)$-bundle $\mathbf{P}$. Though the Levi-Civita connection is not always a connection on $\mathbf{P}$, we can project it onto $\mathbf{P}$ and obtain a connection $\omega$ on $\mathbf{P}$. In other words, $\omega$ is the $\mathfrak{sp}(1)\oplus \mathfrak{sp}(n)$-part of the Levi-Civita connection (cf. \cite{H2}). Note that $M$ is a quaternionic K\"ahler manifold if the torsion tensor of $\omega$ is zero. In the same manner as \eqref{eqn:gradient}, we construct first order differential operators $\{D_{N,\mu}\}_{N,\nu}$ on $\mathbf{S}_{k,\rho}$ with respect to $\omega$. We also call them gradients. 

Let $(M,g)$ be a quaternionic K\"ahler manifold. A conformal deformation $g':=e^{2\sigma(x)}g$ of the Riemannian metric $g$ gives an almost quaternionic Hermitian manifold $(M,g')$. Then we have two gradients, $D_{N,\nu}$ on $(M,g)$ and $D_{N,\nu}'$ on $(M,g')$. The next proposition is the reason why $w_{\nu}$ and $W_N$ are called conformal weights. 
\begin{proposition}\label{proposition:5-1}
The gradient $D_{N,\nu}$ on $(M,g)$ is related to $D_{N,\nu}'$ on $(M,g')=(M,e^{2\sigma}g)$ covariantly,
\begin{equation}
D_{N,\nu}'=\exp \left(\left(-\frac{w_{\nu}}{2}-\frac{W_N}{2n}-1\right)\sigma(x)\right) \circ D_{N,\nu}\circ \exp\left(\left(\frac{w_{\nu}}{2}+\frac{W_N}{2n}\right)\sigma(x)\right). \nonumber
\end{equation}
\end{proposition}
\begin{proof}
In \cite{H2}, the author showed a conformal covariance of gradients on a hyper-K\"ahler manifold. In the same way, we can prove the above conformal covariance of quaternionic K\"ahlerian gradients. \end{proof}

\section{Bochner-Weitzenb\"ock formulas}\label{sec:6}
We consider a vector bundle $\mathbf{S}_{k,\rho}$ on a quaternionic K\"ahler manifold $M$. Set $B_{N,\nu}:=(D_{N,\nu})^{\ast}D_{N,\nu}$, and we know that the second order operator has the following  expression,
\begin{equation}
B_{N,\nu}=(D_{N,\nu})^{\ast}D_{N,\nu}=-\sum_{A,B,\alpha,\beta}\sgn(A\alpha)p_{N,\nu}(v_{B,\beta})^{\ast}p_{N,\nu}(v_{A,\alpha})\nabla^2_{v_{B,\beta},v_{-A,-\alpha}},\nonumber
\end{equation}
where $\nabla^2_{X,Y}$ is defined to be $\nabla_X\nabla_Y-\nabla_{\nabla_X Y}$ for vector fields $X$ and $Y$. It follows from \eqref{eqn:4-4} of $q=0$ that 
\begin{equation}
\sum_{N,\nu}B_{N,\nu}=\nabla^{\ast}\nabla=-\sum_{A,\alpha} \sgn(A\alpha) \nabla^2_{v_{A,\alpha},v_{-A,-\alpha}}, \nonumber
\end{equation}
where $\nabla^{\ast}\nabla$ is the connection Laplacian on $\mathbf{S}_{k,\rho}$. Thus a linear combination of $\{B_{N,\nu}\}_{N,\nu}$ has second order in general. But, some appropriate combinations are curvature endomorphisms which are zeroth order operators, that is,
\[
\sum_{N,\nu} a_{N,\nu}B_{N,\nu}=\textrm{(curvature endomorphism)}.
\]
We call such equations \textit{(optimal) Bochner-Weitzenb\"ock formulas} (see \cite{Br1}). 

Let us give Bochner-Weitzenb\"ock formulas for quaternionic K\"ahlerian gradients. We define the curvature $R_{k,\rho}$ of $\nabla$ on $\mathbf{S}_{k,\rho}$ by
\begin{equation}
R_{k,\rho}(X,Y):=\nabla^2_{X,Y}-\nabla^2_{Y,X}. \nonumber
\end{equation}
Then, from \eqref{relation-1} and \eqref{relation-2}, we have
\begin{equation}
\begin{split}
& \sum_{N,\nu}\left\{(1-(-1)^q)\hat{w}_{\nu}^q+\frac{1-(-1)^q}{2}\hat{w}_{\nu}^{q-1}+\sum_{p=0}^{q-1}(-1)^p \pi_{\rho}(\hat{c}_{q-1-p})\hat{w}_{\nu}^p\right\}B_{N,\nu} \\
&=-\sum_{A,\alpha,\beta} \sgn(A\alpha)\id \hotimes \pi_{\rho}(\hat{x}_{-\alpha\beta}^q)R_{k,\rho}(v_{A,\alpha},v_{-A,\beta}),
 \end{split}\label{bochner-0-1}
\end{equation}
and 
\begin{equation}
\begin{split}
 &\sum_{N,\nu}W_N\left\{(1+(-1)^q)\hat{w}_{\nu}^q-\frac{1-(-1)^q}{2}\hat{w}_{\nu}^{q-1}-\sum_{p=0}^{q-1}(-1)^p \pi_{\rho}(\hat{c}_{q-1-p})\hat{w}_{\nu}^p\right\}B_{N,\nu} \\
&=-\sum_{A,B,\alpha,\beta} \sgn(A\alpha)\pi_k(y_{-AB})\hotimes \pi_{\rho}(\hat{x}_{-\alpha\beta}^q)R_{k,\rho}(v_{A,\alpha},v_{B,\beta}).
\end{split}\label{bochner-0-2}
\end{equation}

To obtain more simple formulas, we should calculate curvature endomorphisms in the above equations. We denote by $R$ the Riemannian curvature tensor on the tangent bundle $T(M)$. In \cite{KSW}, W. Kramer, U. Semmelmann and G. Weingart gave a formula of $R$ as follows. 
We define three $\End(\mathbf{H}\hotimes \mathbf{E})$-valued $2$-forms by
\begin{gather}
R^H(h\hotimes \e,h'\hotimes \e'):=\sigma_E(\e,\e')(h\odot h'\otimes \id_E),\nonumber\\
R^E(h\hotimes \e,h'\hotimes \e'):=\sigma_H(h,h')(\id_H\otimes \e\odot \e'),\nonumber\\
R^{hyper}(h\hotimes \e,h'\hotimes \e'):=\sigma_H(h,h')(\id_E\hotimes \mathfrak{R}(\e,\e')).\nonumber
\end{gather}
Here, $\mathfrak{R}$ is the $S^4(\mathbf{E})$-part of the curvature $R$. In other words, $\sigma_E(\mathfrak{R}(\e^1,\e^2)\e^3,\e^4)$ is symmetric for $\e^1,\e^2,\e^3$ and $\e^4$. Note that $\mathfrak{R}(\e,\e')$ is expressed by
\[
\mathfrak{R}(\e,\e')=1/2\sum_{\delta,\gamma} \sgn(\delta)\sigma_E(\mathfrak{R}(\e,\e')\e_{\gamma},\e_{-\delta})x_{\delta\gamma}.
\]
Then the Riemannian curvature tensor $R$ is 
\[
R=-\frac{\kappa}{8n(n+2)}\left(R^H+R^E\right)+R^{hyper}.
\]
Here, $\kappa$ is the scalar curvature of $R$. Since a quaternionic K\"ahler manifold is an Einstein manifold, the scalar curvature $\kappa$ is constant. Note that $R^{hyper}$ is zero on the quaternionic projective space $\mathbb{H}P^n$ (see \cite{S1}). 

The covariant derivative $\nabla$ on $\mathbf{S}_{k,\rho}$ is defined from the Levi-Civita connection and hence the curvature $R_{k,\rho}$ is $\pi_{k,\rho}(R)$. The curvature endomorphisms on the right sides of \eqref{bochner-0-1} and \eqref{bochner-0-2} are rewritten as follows.
\begin{lemma}
We set
\begin{equation}
\hat{\mathfrak{R}}^q_{\rho}:=\sum_{\alpha,\beta}\sgn(\alpha)\id\hotimes \pi_{\rho}(\hat{x}_{\alpha\beta}^q\mathfrak{R}(\e_{-\alpha},\e_{\beta})). \nonumber
\end{equation}
Then we have 
\begin{equation}
\begin{split}
 &-\sum_{A,\alpha,\beta} \sgn(A\alpha)\id \hotimes \pi_{\rho}(\hat{x}_{-\alpha\beta}^q)R_{k,\rho}(v_{A,\alpha},v_{-A,\beta})\\
=&\frac{\kappa}{4n(n+2)}\pi_{\rho}\left(\hat{c}_{q+1}+\frac{2n+1}{2}\hat{c}_q\right)+2\hat{\mathfrak{R}}^q_{\rho},
\end{split}\nonumber
\end{equation}
and 
\begin{equation}
\begin{split}
 &-\sum_{A,B,\alpha,\beta} \sgn(A\alpha)\pi_k(y_{-AB})\hotimes \pi_{\rho}(\hat{x}_{-\alpha\beta}^q)R_{k,\rho}(v_{A,\alpha},v_{B,\beta})\\
=&\frac{\kappa}{8n(n+2)}\pi_k(C_2)\pi_{\rho}(\hat{c}_q)=\frac{k(k+2)\kappa}{4n(n+2)}\pi_{\rho}(\hat{c}_q).
\end{split}\nonumber
\end{equation}
\end{lemma}

We detect the number of independent Bochner-Weitzenb\"ock formulas. We assume that there are $\mathcal{N}$ irreducible components in $\mathbf{S}_{k,\rho}\otimes (\mathbf{H}\hotimes \mathbf{E})$, that is,
\begin{equation}
\mathcal{N}:=\#\{(k+N,\rho+\mu_{\nu})\;|\; \textrm{both $k+N$ and $\rho+\mu_{\nu}$ are dominant integral}\}.\nonumber
\end{equation}
Then we have $\mathcal{N}$ gradients on $\Gamma(M,\mathbf{S}_{k,\rho})$. From a similar discussion to the one in \cite{H4}, we can show that there are at least $[\mathcal{N}/2]$ independent Bochner-Weitzenb\"ock formulas.
\begin{theorem}\label{theorem:6-2}
We assume that the number of non-zero gradients is $\mathcal{N}$. The operators $\{B_{N,\nu}=(D_{N,\nu})^{\ast}D_{N\nu}\}_{N,\nu}$ on $\Gamma(M, \mathbf{S}_{k,\rho})$ satisfy
\begin{equation}
\sum_{N,\nu}B_{N,\nu}=\nabla^{\ast}\nabla. \nonumber
\end{equation}
Furthermore, when the highest weight $k$ for $\Sp(1)$ is not zero, we have the following $[\mathcal{N}/2]$ independent Bochner-Weitzenb\"ock formulas:
\begin{enumerate} 
\item For $q=1,2,\dots ,[\mathcal{N}/4]$,  
\begin{equation}
\begin{split}
& \sum_{N,\nu}\left\{\sum_{p=0}^{2q-1}(-1)^p \pi_{\rho}(\hat{c}_{2q-1-p})\hat{w}_{\nu}^p\right\}B_{N,\nu} \\
=&\frac{\kappa}{4n(n+2)}\pi_{\rho}\left(\hat{c}_{2q+1}+\frac{2n+1}{2}\hat{c}_{2q}\right)+2\hat{\mathfrak{R}}^{2q}_{\rho}.
 \end{split}\label{bochner-1}
\end{equation} 
\item For $q=0,1,\dots ,[\mathcal{N}/4-1/2]$.
\begin{equation}
\begin{split}
 \sum_{N,\nu}W_N\left\{2\hat{w}_{\nu}^{2q}-\sum_{p=0}^{2q-1}(-1)^p \pi_{\rho}(\hat{c}_{2q-1-p})\hat{w}_{\nu}^p\right\}B_{N,\nu}
=\frac{k(k+2)\kappa}{4n(n+2)}\pi_{\rho}(\hat{c}_{2q}).
\end{split} \label{bochner-2}
\end{equation}
\end{enumerate}
When $k$ is zero, the equations \eqref{bochner-1} for $q=1,\dots,[\mathcal{N}/2]$ constitute $[\mathcal{N}/2]$ independent Bochner-Weitzenb\"ock formulas. 
\end{theorem}
It is useful to write down the first few formulas. From \eqref{eqn:ubw}--\eqref{eqn:casimir-1}, \eqref{eqn:c-1} and \eqref{eqn:c-2}, we rewrite \eqref{bochner-1} of $q=1,2$ as
\begin{gather}
 \sum_{N,\nu}w_{\nu}B_{N,\nu}
=\frac{\kappa}{8n(n+2)}\pi_{\rho}(c_2)+\sum_{\alpha,\beta}\sgn(\alpha)\id\hotimes \pi_{\rho}(x_{\alpha\beta}\mathfrak{R}(\e_{-\alpha},\e_{\beta})), \label{eqn:bw-1} \\
\begin{split}
 &\sum_{N,\nu}\left\{\pi_{\rho}(c_2)/2+(n+1)(2n+1)w_{\nu}-(2n+1)w_{\nu}^2+w_{\nu}^3\right\}B_{N,\nu}\\
&=\frac{\kappa}{8n(n+2)}\pi_{\rho}(c_4)+\sum_{\alpha,\beta}\sgn(\alpha)\id\hotimes \pi_{\rho}(x_{\alpha\beta}^3\mathfrak{R}(\e_{-\alpha},\e_{\beta})),
\end{split}\label{eqn:bw-2}
\end{gather}
and \eqref{bochner-2} of $q=0,1,2$ as
\begin{gather}
\sum_{N,\nu}W_N B_{N,\nu}=\frac{k(k+2)\kappa}{4(n+2)}, \label{eqn:bw-3}\\
\sum_{N,\nu}2W_N\left(w_{\nu}^2-(n+1)w_{\nu}\right)B_{N,\nu}=\frac{k(k+2)\kappa}{4n(n+2)}\pi_{\rho}(c_2), \label{eqn:bw-4} \\
\begin{split}
\sum_{N,\nu}&W_N\left\{ 2w_{\nu}(w_{\nu}-n-1)(w_{\nu}^2-(2n+1)w_{\nu}+2n+1)\right.\\
&+\left.(n+w_{\nu}) \pi_{\rho}(c_2)\right\}B_{N,\nu}=\frac{k(k+2)\kappa}{4n(n+2)}\pi_{\rho}(c_4).
\end{split}\label{eqn:bw-5}
\end{gather}

To apply Bochner-Weitzenb\"ock formulas to differential geometry, we compare the curvature endomorphisms for the Riemannian case and the ones for the quaternionic K\"ahler case. In \cite{H1} and \cite{H4}, the author discussed the Bochner-Weitzenb\"ock formulas and associated curvature endomorphisms on a Riemannian or spin manifold. Let $(M,g)$ be an oriented $m$-dimensional Riemannian manifold. We take an irreducible vector bundle $\mathbf{S}_{\lambda}$ associated to the orthonormal frame bundle of $M$, where $\lambda$ is a dominant integral weight with respect to $\mathrm{SO}(m)$. Then we define a curvature endomorphism $R_{\lambda}^1$ on $\mathbf{S}_{\lambda}$ by
\begin{equation}
R^1_{\lambda}=\sum_{1\le i,j\le m}\pi_{\lambda}(e_i\wedge e_j)R_{\lambda}(e_i,e_j), \nonumber
\end{equation}
where $\{e_i\}_{1\le i\le m}$ is a local orthonormal frame of $M$ and $R_{\lambda}$ is the curvature on $\mathbf{S}_{\lambda}$. When $M$ has a spin structure, we can also define a curvature endomorphism $R^1_{\lambda}$ on vector bundle $\mathbf{S}_{\lambda}$ associated to the bundle of spin frames. 

If $M$ has a quaternionic K\"ahler structure, we can decompose $\mathbf{S}_{\lambda}$ into irreducible bundles with respect to $\Sp(1)\Sp(n)\subset \mathrm{SO}(4n)$, $
\mathbf{S}_{\lambda}=\bigoplus_{k,\rho}\mathbf{S}_{k,\rho}$. 
Then it is  easily to show that the restriction of $R_{\lambda}^1$ onto $\mathbf{S}_{k,\rho}$ is given by 
\begin{equation}
R_{k,\rho}^1:=\frac{\kappa}{8n(n+2)}\left(\pi_k(C_2)+\pi_{\rho}(c_2)\right)+\sum_{\alpha,\beta}\sgn(\alpha)\id\hotimes \pi_{\rho}(x_{\alpha\beta}\mathfrak{R}(\e_{-\alpha},\e_{\beta})).
  \nonumber
\end{equation}
From \eqref{eqn:bw-1} and \eqref{eqn:bw-3}, we have a quaternionic K\"ahlerian version of Gauduchon's formula in \cite{G},
\begin{equation}
R_{k,\rho}^1=R_{\lambda}^1|_{\mathbf{S}_{k,\rho}}=\sum_{N,\nu}( w_{\nu}+W_N/n)B_{N,\nu}. \nonumber
\end{equation}
\begin{remark}
The curvature endomorphism $R_{k,\rho}^1$ corresponds to $4q(R)$ in \cite{SW}. 
\end{remark}
\begin{example}
We consider the bundle $\mathbf{S}_{\lambda}=\Lambda^p(M)$. Then we have a Bochner-Weitzenb\"ock formula for the Laplace operator (see \cite{H4}), 
\[
dd^{\ast}+d^{\ast}d=\nabla^{\ast}\nabla+R_{\lambda}^1/2.
\]
Restricting the above equation onto an irreducible bundle $\mathbf{S}_{k,\rho}$, we have
\[
dd^{\ast}+d^{\ast}d=\nabla^{\ast}\nabla+R_{k,\rho}^1/2=\sum_{N,\nu}\left(1+\frac{w_{\nu}}{2}+\frac{W_N}{2n}\right)B_{N,\nu} .
\]
\end{example}
\begin{example}
We consider the quaternionic projective space $\mathbb{H}P^n$ with $\kappa=2n$. Since the curvature $R^{hyper}=0$, we obtain
\[
dd^{\ast}+d^{\ast}d=\nabla^{\ast}\nabla+R^1_{k,\rho}/2=\nabla^{\ast}\nabla+\frac{1}{8(n+2)}(2k(k+2)+\pi_{\rho}(c_2)).
\]
\end{example}
\begin{example}
Let $M$ be a quaternionic K\"ahler spin manifold and $\mathbf{S}_{\lambda}$ be the spinor bundle. The Dirac operator satisfies 
\[
D^2=\nabla^{\ast}\nabla+R^1_{\lambda}=\nabla^{\ast}\nabla+\kappa/4.
\]
With quaternionic K\"ahler structure of $M$, we decompose the spinor bundle as
\[
\mathbf{S}_{\lambda}=\bigoplus_{1\le k\le n} \mathbf{S}_{k,(1_{n-k})}=\bigoplus_{1\le k\le n} S^k(\mathbf{H})\hotimes \Lambda^{n-k}_0(\mathbf{E}).
\]
The restriction of $D^2$ to $\mathbf{S}_{k,(1_{n-k})}$ is 
\[
D^2=\nabla^{\ast}\nabla+\kappa/4=\nabla^{\ast}\nabla+R^1_{k,(1_{n-k})}=\sum_{N,\nu}(1+w_{\nu}+W_N/n)B_{N,\nu} .
\]
\end{example}

\section{Clifford algebras and curvature endomorphisms on differential forms}\label{sec:7}
We define the curvature endomorphisms $\{\mathfrak{R}_{\rho}^q\}_{q\ge 0}$ by 
\begin{equation}
\begin{split}
\mathfrak{R}_{\rho}^q:=&\sum_{\alpha,\beta}\sgn(\alpha)\id\hotimes \pi_{\rho}(x_{\alpha\beta}^q\mathfrak{R}(\e_{-\alpha},\e_{\beta}))\\
 =&1/2\sum_{\alpha,\beta,\delta,\gamma}\sgn(\alpha\delta)\sigma_E(\mathfrak{R}(\e_{-\alpha},\e_{\beta})\e_{\gamma},\e_{-\delta})\id \hotimes \pi_{\rho}(x_{\alpha\beta}^q x_{\delta\gamma}).
\end{split} \nonumber
\end{equation}
We shall investigate the above endomorphisms on $\Lambda^p(M)$ more precisely. It is known that each irreducible summand in $\Lambda^p(M)$ is the form of $\mathbf{S}_{k,(2_b,1_{a-b})}=S^k(\mathbf{H})\hotimes \Lambda^{a,b}_0(\mathbf{E})$ for $0\le k\le 2n-a-b$ and $0\le b \le a\le n$ (cf. \cite{SW}). Since the curvature endomorphism $\mathfrak{R}^q_{\rho}$ is independent of the $\Sp(1)$-part, we may consider only $\mathbf{S}_{0,(2_b,1_{a-b})}=\Lambda^{a,b}_0(\mathbf{E})$. The main result in this section is the following. 
\begin{proposition}\label{proposition:7-1}
On $\Lambda^{a,b}_0(\mathbf{E})$, we have 
\begin{equation}
\mathfrak{R}^3_{(2_b,1_{a-b})}=(2n^2+7n+7-\pi_{(2_b,1_{a-b})}(c_2)/4)\mathfrak{R}^1_{(2_b,1_{a-b})}.  \label{eqn:curvature-sp}
\end{equation}
Furthermore, we have $\mathfrak{R}^1_{(1_{a})}=0$ on $\Lambda^a_0(\mathbf{E})=\Lambda^{a,0}_0(\mathbf{E})$. 
\end{proposition}
Our method is to make use  of the Clifford algebra $\Cl(E)$ associated to $(E,\sigma_E,J_E)$. Let $\{\e_{\alpha}\}_{\alpha}$ be a symplectic unitary basis of $E$. The Clifford algebra $\Cl(E)$ is an associated algebra over $\mathbb{C}$ generated by $\{\e_{\alpha},\e_{\alpha}^{\dagger}\}_{\alpha}\cup \{1\}$ with relations
\[
\e_{\alpha}\e_{\beta}+\e_{\beta}\e_{\alpha}=0,\quad \e_{\alpha}^{\dagger}\e_{\beta}^{\dagger}+\e_{\beta}^{\dagger}\e_{\alpha}^{\dagger}=0,\quad \e_{\alpha}\e_{\beta}^{\dagger}+\e_{\beta}^{\dagger}\e_{\alpha}=\delta_{\alpha,\beta}.
\]
This algebra acts on $\bigoplus_{a=0}^n \Lambda^a(E)$ by 
\[
\e_{\alpha}^{\dagger}\cdot:=i(\e_{\alpha}),  \quad \e_{\alpha}\cdot:=\e_{\alpha}\wedge,
\]
where $i(\e_{\alpha})$ is the interior product and $\e_{\alpha}\wedge$ is the exterior product. Then the representation $(\pi_{a}, \Lambda^a(E))$ of $\sypc$ is realized by 
\[
\pi_a(x_{\alpha\beta})=\e_{\alpha}\e_{\beta}^{\dagger}-\sgn(\alpha\beta)\e_{-\beta}\e_{-\alpha}^{\dagger}
\]
for $x_{\alpha\beta}$ in $\sypc$. In other words, the Lie algebra $\sypc$ is embedded in $\Cl(E)$ by $x_{\alpha\beta}\mapsto \e_{\alpha}\e_{\beta}^{\dagger}-\sgn(\alpha\beta)\e_{-\beta}\e_{-\alpha}^{\dagger}$. 
\begin{remark}
There are $\sypc$-invariant elements in $\Cl(E)$, 
\begin{gather}
\mathbf{N}:=\sum_{\alpha}\e_{\alpha}\e_{\alpha}^{\dagger},\quad \mathbf{N}^{\dagger}:=\sum_{\beta}\e_{\alpha}^{\dagger}\e_{\alpha}, \nonumber\\
\sigma_E:=\sum \sgn(\alpha)\e_{\alpha}\e_{-\alpha},\quad \sigma^{\dagger}_E:=\sum \sgn(\alpha)\e_{\alpha}^{\dagger}\e_{-\alpha}^{\dagger}.\nonumber
\end{gather}
The operator $\mathbf{N}$ is the so-called number operator acting on $\Lambda^a(E)$ by constant $a$, and $\mathbf{N}^{\dagger}=2n-\mathbf{N}$. The operator $\sigma_E$ is used for decomposing $\Lambda^a(E)$. More precisely, we have $[\sigma_E,\sigma^{\dagger}_E]=4n-4a$ on $\Lambda^a(E)$, and $\Lambda^a(E)=\bigoplus_{p=0}^{[a/2]}\sigma_E^p\Lambda^{a-2p}_0(E)$. We use these invariant operators implicitly in the proof below. 
\end{remark}
\begin{proof}[Proof of Proposition \ref{proposition:7-1}]
Though the second claim that $\mathfrak{R}^1_{(1_{a})}$ is zero has been already shown in \cite{KSW}, we give a proof of it as a good exercise before proving the first claim. 

It follows from the Clifford relations and the symmetry of $\sigma_E(\mathfrak{R}(\cdot,\cdot)\cdot,\cdot)$ that 
\begin{equation}
\begin{split}
 &\sum_{\alpha}\sgn(\alpha)\pi_a(x_{\alpha\beta}\mathfrak{R}(\e_{-\alpha},\e_{\beta})) \\
=&1/2\sum_{\alpha,\beta,\delta,\gamma}\sgn(\alpha\delta)\sigma_E(\mathfrak{R}(\e_{-\alpha},\e_{\beta})\e_{\gamma},\e_{-\delta}) \pi_{a}(x_{\alpha\beta} x_{\delta\gamma}) \\
=&1/2\sum_{\alpha,\beta,\delta,\gamma}\sgn(\alpha\delta)\sigma_E(\mathfrak{R}(\e_{-\alpha},\e_{\beta})\e_{\gamma},\e_{-\delta})\\
&\qquad (\e_{\alpha}\e_{\beta}^{\dagger}-\sgn(\alpha\beta)\e_{-\beta}\e_{-\alpha}^{\dagger})(\e_{\delta}\e_{\gamma}^{\dagger}-\sgn(\delta\gamma)\e_{-\gamma}\e_{-\delta}^{\dagger})\\
=&1/2\sum_{\alpha,\beta,\delta,\gamma}\sgn(\alpha\delta)\sigma_E(\mathfrak{R}(\e_{-\alpha},\e_{\beta})\e_{\gamma},\e_{-\delta}) \\
&\{\sgn(\beta\delta)(\delta_{\beta,\delta}\e_{\alpha}\e_{\gamma}^{\dagger}+\delta_{-\alpha,\delta}\e_{-\beta}\e_{\gamma}^{\dagger}+\delta_{\beta,-\gamma}\e_{\alpha}\e_{-\delta}^{\dagger}+\delta_{\alpha,\gamma}\e_{-\beta}\e_{-\delta}^{\dagger})\\
&\quad -\e_{\alpha}\e_{\delta}\e_{\beta}^{\dagger}\e_{\gamma}^{\dagger}+\sgn(\alpha\beta)\e_{-\beta}\e_{\delta}\e_{-\alpha}^{\dagger}\e_{\gamma}^{\dagger}\\
&\qquad +\sgn(\delta\gamma)\e_{\alpha}\e_{-\gamma}\e_{\beta}^{\dagger}\e_{-\delta}^{\dagger}-\sgn(\alpha\beta\delta\gamma)\e_{-\beta}\e_{-\gamma}\e_{-\alpha}^{\dagger}\e_{-\delta}^{\dagger}\}\\
=&0.
\end{split}\nonumber
\end{equation}
Thus we have proved the second claim, $\mathfrak{R}^1_{(1_a)}=0$. To prove the first claim, we consider the tensor representation $(\pi_a\otimes \pi_b,\Lambda^a(E)\otimes \Lambda^b(E))$. From a tedious calculation, we show 
\[
\begin{split}
 &(\pi_a\otimes \pi_b)\left(\sum_{\alpha,\beta}\sgn(\alpha)x_{\alpha\beta}^3\mathfrak{R}(\e_{-\alpha},\e_{\beta})\right)\\
=&2(2n^2+7n+7-(\pi_a\otimes \pi_b)(c_2/4))\left(\sum_{\alpha,\beta}\sgn(\alpha)\pi_a(\mathfrak{R}(\e_{-\alpha},\e_{\beta}))\otimes \pi_b(x_{\alpha\beta})\right).
\end{split}
\]
On the other hand, we have
\[
\begin{split}
 (\pi_a\otimes \pi_b)\left(\sum_{\alpha,\beta}\sgn(\alpha)x_{\alpha\beta}\mathfrak{R}(\e_{-\alpha},\e_{\beta})\right)=2\sum_{\alpha,\beta}\sgn(\alpha)\pi_a(\mathfrak{R}(\e_{-\alpha},\e_{\beta}))\otimes \pi_b(x_{\alpha\beta}).
\end{split}
\]
Then we conclude that
\begin{multline}
(\pi_a\otimes \pi_b)\left(\sum_{\alpha,\beta}\sgn(\alpha)x_{\alpha\beta}^3\mathfrak{R}(\e_{-\alpha},\e_{\beta})\right)\\
=(2n^2+7n+7-(\pi_a\otimes \pi_b)(c_2/4))
(\pi_a\otimes \pi_b)\left(\sum_{\alpha,\beta}\sgn(\alpha)x_{\alpha\beta}\mathfrak{R}(\e_{-\alpha},\e_{\beta})\right).\nonumber
\end{multline}
Restricting this equation onto $\Lambda^{a,b}_0(\mathbf{E})$ in $\Lambda^a(\mathbf{E})\otimes \Lambda^b(\mathbf{E})$, we have proved the first claim. 
\end{proof}
By using \eqref{eqn:curvature-sp}, we eliminate the curvature endomorphism $\mathfrak{R}^1_{(2_b,1_{a-b})}$ from \eqref{eqn:bw-1} and \eqref{eqn:bw-2}.
\begin{corollary}
On $S^k(\mathbf{H})\hotimes \Lambda^{a,b}_0(\mathbf{E})$, there is a Bochner-Weitzenb\"ock formula depending only on the scalar curvature $\kappa$, 
\begin{equation}
\begin{split}
 &\sum_{N,\nu}(w_{\nu}+2)\left(\pi_{(2_b,1_{a-b})}(c_2)+4w_{\nu}^2-8nw_{\nu}-12w_{\nu}\right)B_{N,\nu}\\
=&\frac{\kappa}{8n(n+2)}\pi_{(2_b,1_{a-b})}\left(-4(2n^2+7n+7)c_2+c_2^2+4c_4\right).
\end{split}\label{eqn:bw-6}
\end{equation}
This formula is linear independent of \eqref{eqn:bw-3}--\eqref{eqn:bw-5}. 
\end{corollary}

\section{Eigenvalue estimates}\label{sec:8}
We shall apply our Bochner-Weitzenb\"ock formulas to eigenvalue estimate on $\mathbf{S}_{k,(2_b,1_{a-b})}=S^k(\mathbf{H})\hotimes \Lambda^{a,b}_0(\mathbf{E})$ for $0\le k\le 2n-a-b$ and $0\le b\le a\le n$. This bundle is an irreducible summand of the bundle of differential forms \cite{SW}. We will discuss four cases in turn: (1) $a=b=0$, (2) $a>b=0$, (3) $a=b>0$, (4) $a>b>0$. In conclusion, the eigenvalue estimate of the Laplace operator $dd^{\ast}+d^{\ast}d$ restricted on $\mathbf{S}_{k,(2_b,1_{a-b})}$ is given in Theorem \ref{theorem:8-2}. 

\subsection{Estimates on $S^k(\mathbf{H})$}\label{sec:8-0}
We consider the bundle $\mathbf{S}_{k,(0_n)}=S^k(\mathbf{H})$. There are two gradients, $D_{1,1}$ and $D_{-1,1}$. These operators satisfy
\[
B_{1,1}+B_{-1,1}=\nabla^{\ast}\nabla,\quad -k B_{1,1}+(k+2)B_{-1,1}=\frac{k(k+2)}{4(n+2)}\kappa,
\]
where $B_{\pm 1,1}=(D_{\pm 1,1})^{\ast}D_{\pm 1,1}$. Then we have 
\[
\nabla^{\ast}\nabla=\frac{2(k+1)}{k+2}B_{1,1}+\frac{k}{4(n+2)}\kappa =\frac{2(k+1)}{k}B_{-1,1}-\frac{k+2}{4(n+2)}\kappa.
\]
We think of $S^k(\mathbf{H})$ as an irreducible summand of the bundle of differential forms. The restricted Laplacian $dd^{\ast}+d^{\ast}d$ on $S^k(\mathbf{H})$ is 
\[
\nabla^{\ast}\nabla+R^1_{k,(0_n)}/2=\nabla^{\ast}\nabla+\frac{k(k+2)}{8n(n+2)}\kappa.
\]
Then we have 
\[
\begin{split}
dd^{\ast}+d^{\ast}d &=\frac{2(k+1)}{k+2}B_{1,1}+\frac{k(2n+k+2)}{8n(n+2)}
   \kappa\\
 &=\frac{2(k+1)}{k}B_{-1,1}-\frac{(k+2)(2n-k)}{8n(n+2)}\kappa.
\end{split}
\]
This equation leads to eigenvalue estimates of $dd^{\ast}+d^{\ast}d$ (cf. \cite{SW2}).
\begin{proposition}
We consider a compact quaternionic K\"ahler manifold with non-zero scalar curvature. A lower bound of the eigenvalues of $dd^{\ast}+d^{\ast}d$ on $S^k(\mathbf{H})$ for non-negative integer $k$ is given as follows.
\[
\begin{cases}
\displaystyle \frac{k(2n+k+2)}{8n(n+2)}\kappa & \textrm{for $\kappa>0$},\\
\displaystyle -\frac{(k+2)(2n-k)}{8n(n+2)}\kappa & \textrm{for $\kappa<0$}.
\end{cases}
\]
\end{proposition}

\subsection{Estimates on $S^k(\mathbf{H})\hotimes \Lambda^a_0(\mathbf{E})$}\label{sec:8-1}
We consider $\mathbf{S}_{k,(1_{a})}=S^k(\mathbf{H})\hotimes \Lambda^{a}_0(\mathbf{E})$ for $a\neq 0$. Then we have six gradients on this vector bundle, $D_{1,1}$, $D_{1,a+1}$, $D_{1,-a}$, $D_{-1,1}$, $D_{-1,a+1}$ and $D_{-1,-a}$, where we set $D_{\pm 1,a+1}:=0$ in the case of $a=n$. 

The sum of $B_{N,\nu}=(D_{N,\nu})^{\ast}D_{N,\nu}$ is the connection Laplacian,
\begin{equation}
B_{1,1}+B_{1,a+1}+B_{1,-a}+B_{-1,1}+B_{-1,a+1}+B_{-1,-a}=\nabla^{\ast}\nabla.
\label{eqn:bw-a-1}
\end{equation}
Moreover, it follows from \eqref{eqn:bw-1}, \eqref{eqn:bw-3} and \eqref{eqn:bw-4} that there are three independent Bochner-Weitzenb\"ock formulas, 
\begin{gather}
\begin{split}
&-B_{1,1}+ a  B_{1,a+1}+(2n-a+2)B_{1,-a}\\
 &\quad -B_{-1,1}+a B_{-1,a+1}+(2n-a+2)B_{-1,-a}=\frac{a(2n-a+2)}{4n(n+2)}\kappa,\end{split} \nonumber \\
\begin{split}
&-k(B_{1,1}+B_{1,a+1}+B_{1,-a})\\
&\quad+(k+2)(B_{-1,1}+B_{-1,a+1}+B_{-1,-a}) =\frac{k(k+2)}{4(n+2)}\kappa,
\end{split}\label{eqn:bw-a-2}\\
\begin{split}
&-k(n+2)B_{1,1}+ ka(n-a+1)B_{1,a+1}-k(2n-a+2)(n-a+1)B_{1,-a}\\
&+(k+2)(n+2)B_{-1,1}-(k+2)a(n-a+1)B_{-1,a+1}\\
&+(k+2)(2n-a+2)(n-a+1)B_{-1,-a}=\frac{k(k+2)a(2n-a+2)}{4n(n+2)}\kappa.
\end{split}\nonumber
\end{gather}
From these formulas, we change the connection Laplacian $\nabla^{\ast}\nabla$ into the form of
\begin{equation}
\sum_{N,\nu} c_{N,\nu}B_{N,\nu}+c\kappa, \quad \textrm{ for $c_{N,\nu}\ge 0$}. \label{eqn:laplacian}
\end{equation}
Since $B_{N,\nu}$ is non-negative operator on a compact quaternionic K\"ahler manifold, the eigenvalues of $\nabla^{\ast}\nabla$ have a lower bound $c\kappa$. For example, it follows from \eqref{eqn:bw-a-1} and \eqref{eqn:bw-a-2} that 
\begin{equation}
\nabla^{\ast}\nabla=\frac{2(k+1)}{k+2}(B_{1,1}+B_{1,a+1}+B_{1,-a})+\frac{k}{4(n+2)}\kappa. \nonumber
\end{equation}
Then the eigenvalues of $\nabla^{\ast}\nabla$ on a compact positive quaternionic K\"ahler manifold have a lower bound $\frac{k}{4(n+2)}\kappa$. Thus, to get a lower bound of the eigenvalues of $\nabla^{\ast}\nabla$, we should find out a formula of \eqref{eqn:laplacian} such that $c\kappa$ is as great as possible. In fact we can rewrite $\nabla^{\ast}\nabla$ as follows.
\begin{enumerate}
\item In the case that the scalar curvature $\kappa$ is positive,
\begin{enumerate}
\item For $k=0$,
\begin{equation}
\begin{split}
\nabla^{\ast}\nabla=\frac{(2n-a+3)}{(2n-a+2)}B_{1,1}+\frac{2(n-a+1)}{(2n-a+2)}B_{1,a+1}+\frac{a}{4n(n+2)}\kappa.
\end{split}\nonumber
\end{equation}
\item For $k\neq 0$,
\[
\nabla^{\ast}\nabla=\frac{2(k+1)}{k+2}(B_{1,1}+B_{1,a+1}+B_{1,-a})+\frac{k}{4(n+2)}\kappa.
\]
\end{enumerate}
\item In the case that the scalar curvature $\kappa$ is negative,
\begin{enumerate}
\item For $0\le k\le n-a$,
\begin{equation}
\begin{split}
&\nabla^{\ast}\nabla=\frac{2(a+1)(n-k-a)}{(k+2)(n-a)}B_{1,a+1}+\frac{2(k+1)(2n-a+3)}{k+2}B_{1,-a}\\
&+\frac{2(a+1)(n-a+1)}{n-a}B_{-1,a+1}-\frac{2an+kn-a^2-ka+2a+2k}{4n(n+2)}\kappa. \end{split}\nonumber
\end{equation}
\item For $n-a< k\le 2n-a$,
\begin{equation}
\begin{split}
&\nabla^{\ast}\nabla =\frac{2(2n-a+3)(n-a+1)}{n-a+2}B_{1,-a}+\frac{2(k+1)(a+1)}{k}B_{-1,a+1}\\ 
&+\frac{2(2n-a+3)(k+a-n)}{k(n-a+2)}B_{-1,-a}-\frac{-ka-a^2+2n+kn+2an}{4n(n+2)}\kappa.
\end{split}\nonumber
\end{equation}
\end{enumerate}
\end{enumerate}
Then we have an eigenvalue estimate of $\nabla^{\ast}\nabla$ on $\mathbf{S}_{k,(1_a)}$. 
\begin{proposition}\label{proposition:8-1}
We consider the connection Laplacian $\nabla^{\ast}\nabla$ on $\mathbf{S}_{k,(1_a)}=S^k(\mathbf{H})\hotimes \Lambda^a_0(\mathbf{E})$. A lower bound of the eigenvalues of $\nabla^{\ast}\nabla$ is as follows.
\begin{enumerate}
\item On a compact positive quaternionic K\"ahler manifold, 
\[
\begin{cases}
\displaystyle\frac{a}{4n(n+2)}\kappa & \textrm{for $k=0$},\\
\displaystyle\frac{k}{4(n+2)}\kappa & \textrm{for $k\neq 0$}.
\end{cases}
\]
\item On a compact negative quaternionic K\"ahler manifold, 
\[
\begin{cases}
\displaystyle -\frac{2an+kn-a^2-ka+2a+2k}{4n(n+2)}\kappa & 
  \textrm{for $0\le k\le n-a$}, \\
\displaystyle-\frac{-ka-a^2+2n+kn+2an}{4n(n+2)}\kappa & \textrm{for $n-a< k\le 2n-a$}.
\end{cases}
\]
\end{enumerate}
\end{proposition}

\begin{example}[The Dirac operator \cite{KSW}]
We consider the Dirac operator on the spinor bundle $\bigoplus_{k=0}^n\mathbf{S}_{k,(1_{n-k})}$. Because of $D^2=\nabla^{\ast}\nabla+\kappa/4$, the eigenvalues of $D^2$ on a compact positive quaternionic K\"ahler spin manifold have the following lower bound, 
\[
\begin{cases}
\displaystyle\frac{n+3}{4(n+2)}\kappa & \textrm{for $k=0$},\\
\displaystyle\frac{n+k+2}{4(n+2)}\kappa & \textrm{for $0<k\le n$}.
\end{cases}
\]
\end{example}
\begin{example}[The Laplacian]
We think of $\mathbf{S}_{k,(1_a)}$ as an irreducible summand of the bundle of differential forms. The restricted Laplacian $dd^{\ast}+d^{\ast}d$ on $\mathbf{S}_{k,(1_a)}$ is 
\[
\nabla^{\ast}\nabla+R^1_{k,(1_a)}/2=\nabla^{\ast}\nabla+\frac{\kappa}{8n(n+2)}(k(k+2)+a(2n-a+2)).
\]
Then we have eigenvalue estimates of $dd^{\ast}+d^{\ast}d$ on a compact quaternionic K\"ahler manifold with non-zero scalar curvature. 
\begin{enumerate}
\item When $\kappa>0$, a lower bound of the eigenvalues  of $dd^{\ast}+d^{\ast}d$ is 
\[
\begin{cases}
\displaystyle \frac{a(2n-a+4)}{8n(n+2)}\kappa & \textrm{for $k=0$},\\
\displaystyle \frac{(a+k)(2n-a+k+2)}{8n(n+2)}\kappa & \textrm{for $k\neq 0$}.
\end{cases}
\]
\item When $\kappa<0$, a lower bound is
\[
\begin{cases}
\displaystyle -\frac{(a+k)(2n-a-k+2)}{8n(n+2)}\kappa & \textrm{for $0\le k\le n-a$},\\
\displaystyle -\frac{(a+k+2)(2n-a-k)}{8n(n+2)}\kappa &\textrm{for $n-a< k\le 2n-a$}.
\end{cases}
\]
\end{enumerate}
These results give eigenvalue estimates of $dd^{\ast}+d^{\ast}d$ on $\Lambda^1(M)=\mathbf{H}\hotimes \mathbf{E}$ in \cite{AM}, \cite{L}, \cite{SW}. Next we consider the bundle of $2$-forms, 
\[
\Lambda^2(M)=\mathbf{S}_{2,(0_n)}\oplus \mathbf{S}_{2,(1_2,0_{n-2})}\oplus \mathbf{S}_{0,(2,0_{n-1})}.
\]
A lower bound of the eigenvalues of $dd^{\ast}+d^{\ast}d$ on $\mathbf{S}_{2,(1_2,0_{n-2})}$ is 
\[
\begin{cases}
\displaystyle \frac{n+1}{n(n+2)}\kappa & \textrm{for $\kappa>0$}, \\
 \displaystyle -\frac{n-1}{n(n+2)}\kappa & \textrm{for $\kappa<0$ and $n\ge 4$},\\
\displaystyle -\frac{3(n-2)}{2n(n+2)}\kappa & \textrm{for $\kappa<0$ and $n=2,3$}.
\end{cases}
\]
We know that our estimates for $\kappa<0$ are better than the ones in \cite{SW}.
\end{example}
Now, we shall apply Bochner-Weitzenb\"ock formulas to vanishing theorems. We consider the vector bundle $S^{k+1}(\mathbf{H})\hotimes \mathbf{E}$ on a compact quaternionic K\"ahler manifold $M$. It follows from the Penrose transform that the Dolbeault cohomology $H^1(Z,\mathcal{O}(k))$ ($k\ge 0$) on the twistor space $Z$ of $M$ is isomorphic to the solution space of certain linear differential equation on $M$(see \cite{Ho}, \cite{NN}).  We can easily show that the solution space $\mathcal{S}$ is given by
\[
\mathcal{S}=\ker D_{1,2}\cap \ker D_{1,-1}\cap \ker D_{-1,-1}\subset \Gamma(M,S^{k+1}(\mathbf{H})\hotimes \mathbf{E}).
\]
Because of \eqref{eqn:bw-a-2}, a solution $\phi$ in $\mathcal{S}$ satisfies 
\begin{gather}
\|D_{1,1}\phi\|^2=\int_M\langle D_{1,1}\phi,D_{1,1}\phi\rangle \textrm{vol}_g=-\frac{(k+3)(2n+k+2)}{8n(n+2)(k+2)}\kappa\|\phi\|^2,\label{eqn:v-1}\\ 
\|D_{-1,1}\phi\|^2=\frac{k(k+1)}{8(n+2)(k+2)}\kappa\|\phi\|^2,\nonumber\\ 
\|D_{-1,2}\phi\|^2=\frac{(k+1)(n-1)}{8n(n+2)}\kappa\|\phi\|^2. \label{eqn:v-3}
\end{gather}
If the scalar curvature is negative, then the equation \eqref{eqn:v-3} yields $\mathcal{S}=\{0\}$. This vanishing was shown in \cite{Ho}. When the scalar curvature is positive, we also have $\mathcal{S}=\{0\}$ by \eqref{eqn:v-1} (cf. \cite{L0}). 
\begin{proposition}\label{proposition:vanish}
Let $Z$ be the twistor space of a compact quaternionic K\"ahler manifold with non-zero scalar curvature. Then we have $H^1(Z,\mathcal{O}(k))=0$ for non-negative integer $k$. 
\end{proposition}
The author expect that our Bochner-Weitzenb\"ock formulas  will produce vanishing theorems for higher cohomology $H^i(Z,\mathcal{O}(k))$, for example, vanishing theorems in \cite{NN}.

\subsection{Estimates on $S^k(\mathbf{H})\hotimes \Lambda^{a,a}_0(\mathbf{E})$}\label{sec:8-2}
We consider the vector bundle $\mathbf{S}_{k,(2_a)}=S^k(\mathbf{H})\hotimes \Lambda^{a,a}_0(\mathbf{E})$ as an irreducible summand of the bundle of differential forms, where we assume that $a$ is not zero. On this vector bundle, we have six gradients, $D_{\pm 1,1}$, $D_{\pm 1,a+1}$ and $D_{\pm 1,-a}$. The Laplace operator $dd^{\ast}+d^{\ast}d$ restricted to $\mathbf{S}_{k,(2_a)}$ is 
\begin{equation}
\begin{split}
 &\nabla^{\ast}\nabla+R^1_{k,(2_a)}/2=\sum_{N,\nu}\left(1+\frac{w_{\nu}}{2}+\frac{W_N}{2n}\right)B_{N,\nu}\\
=&-\frac{k}{2n}B_{1,1}+\frac{2n+an-k}{2n}B_{1,a+1}+\frac{2n^2+5n-an-k}{2n}B_{1,-a}\\
&+\frac{k+2}{2n}B_{-1,1}+\frac{2n+an+k+2}{2n}B_{-1,a+1}+\frac{2n^2+5n-an+k+2}{2n}B_{-1,-a}.
\end{split}\nonumber
\end{equation}
There are three independent Bochner-Weitzenb\"ock formulas \eqref{eqn:bw-1}, \eqref{eqn:bw-3} and \eqref{eqn:bw-4}. Since \eqref{eqn:bw-1} includes $\mathfrak{R}^1_{k,(2_a)}$, we use only \eqref{eqn:bw-3} and \eqref{eqn:bw-4} to estimate eigenvalues. We rewrite $dd^{\ast}+d^{\ast}d$ to the form of \eqref{eqn:laplacian} such that  $c\kappa$ is greatest. 
\begin{enumerate}
\item In the case of $\kappa>0$, 
\begin{equation}
\begin{split}
 &dd^{\ast}+d^{\ast}d\\
=&\frac{(k+1)(a+2)}{k+2}B_{1,a+1}+\frac{(2n-a+5)(2n-2a+3k+6)}{2(k+2)(n-a+3)}B_{1,-a}\\
&+\frac{(2n-a+5)(2n-2a+3)}{2(n-a+3)}B_{-1,-a}+\frac{k(2n-2a+k+2)}{8n(n+2)}\kappa.
\end{split}\nonumber
\end{equation}
\item In the case of $\kappa<0$, 
\begin{enumerate}
\item For $0\le k\le \frac{2n-2a}{3}$,
\begin{equation}
\begin{split}
 &dd^{\ast}+d^{\ast}d\\
=&\frac{(a+2)(2n-2a-3k)}{2(k+2)(n-a)}
  B_{1,a+1}+\frac{(k+1)(2n-a+5)}{k+2}B_{1,-a} \\
&+\frac{(a+2)(2n-2a+3)}{2(n-a)} B_{-1,a+1}-\frac{k(2n-2a-k+4)}{8n(n+2)}\kappa.
\end{split}\nonumber
\end{equation}
\item For $ \frac{2n-2a}{3}<k\le 2n-2a$, 
\begin{equation}
\begin{split}
&dd^{\ast}+d^{\ast}d\\
=&\frac{(2n-a+5)(2n-2a+3)}{2(n-a+3)}B_{1,-a}+\frac{(k+1)(a+2)}{k}B_{-1,a+1}\\
&+\frac{(2n-a+5)(3k+2a-2n)}{2k(n-a+3)}B_{-1,-a}-\frac{(k+2)(2n-2a-k)}{8n(n+2)}\kappa.
\end{split}\nonumber
\end{equation}
\end{enumerate}
\end{enumerate}
\begin{proposition}
A lower bound of the eigenvalues of $dd^{\ast}+d^{\ast}d$ on $\mathbf{S}_{k,(2_a)}$ is as follows.
\begin{enumerate}
\item On a compact positive quaternionic K\"ahler manifold,
\[
\frac{k(2n-2a+k+2)}{8n(n+2)}\kappa.
\]
\item On a compact negative quaternionic K\"ahler manifold,
\[
\begin{cases}
\displaystyle -\frac{k(2n-2a-k+4)}{8n(n+2)}\kappa & \textrm{for $0\le k\le \frac{2n-2a}{3}$},\\
\displaystyle -\frac{(k+2)(2n-2a-k)}{8n(n+2)}\kappa & \textrm{for $\frac{2n-2a}{3}<k\le 2n-2a$}. 
\end{cases}
\]
\end{enumerate}
\end{proposition}
\subsection{Estimates on $S^k(\mathbf{H})\hotimes \Lambda^{a,b}_0(\mathbf{E})$}
We consider $\mathbf{S}_{k,(2_b,1_{a-b})}=S^k(\mathbf{H})\hotimes \Lambda^{a,b}_0(\mathbf{E})$ for $a>b>0$. We have ten gradients on this bundle,
\[
D_{\pm 1,1},\quad D_{\pm 1,b+1},\quad D_{\pm 1,a+1},\quad D_{\pm 1,-b},\quad D_{\pm 1,-a}.
\]
The Laplacian $dd^{\ast}+dd^{\ast}$ is 
\[
dd^{\ast}+d^{\ast}d=\nabla^{\ast}\nabla+R^1_{k,(2_b,1_{a-b})}/2=\sum_{N,\nu} \left(1+\frac{w_{\nu}}{2}+\frac{W_N}{2n}\right)B_{N,\nu},
\]
and there are four Bochner-Weitzenb\"ock formulas \eqref{eqn:bw-3}--\eqref{eqn:bw-5} and \eqref{eqn:bw-6} depending only on the scalar curvature. By a tedious calculation, we can rewrite $dd^{\ast}+d^{\ast}d$ as follows: \\
On a positive quaternionic K\"ahler manifold, 
\begin{enumerate}
\item For $k=0$, 
\begin{equation}
\begin{split}
 &dd^{\ast}+d^{\ast}d\\
=&\frac{(b+1)(2n-a-b+3)}{2n-a-b+2}B_{1,b+1}+
\frac{2(a+2)(n-a+1)}{(a-b+2)(2n-a-b+2)}B_{1,a+1}\\
&+\frac{(a-b+1)(2n-b+5)}{a-b+2}B_{1,-b}+\frac{(a-b)(2n-a-b+4)}{8n(n+2)}\kappa.
\end{split}\nonumber
\end{equation}
\item For $k\neq 0$, 
\begin{equation}
\begin{split}
&dd^{\ast}+d^{\ast}d\\
=&\frac{2(b+1)(k+1)}{k+2}B_{1,b+1}+\frac{2(a+2)(k+1)}{(k+2)(a-b+2)}B_{1,a+1}\\
&+\frac{2(2n-b+5)}{(a-b+2)(k+2)(2n-a-b+4)(n-b+3)}\times \\
 &(12+6 a - 3 a^2 - 16b - 2a b + a^2 b + 7b^2 - b^3 + 6 k +   2 a k - a^2 k\\
 & - 5 b k + b^2 k + 10 n + 8 a n - a^2 n -12 b n - 2 a b n + 3 b^2 n \\
 &+3 k n + 2a k n - 2 b k n +   2 n^2 + 2 a n^2 - 2 b n^2) B_{1,-b}\\
&+\frac{2(k+1)(2n-a+4)}{(k+2)(2n-a-b+4)}B_{1,-a}\\
&+\frac{2(a-b+1)(2n-b+5)(2n-a-b+3)(n-b+2)}{(a-b+2)(2n-a-b+4)(n-b+3)}
 B_{-1,-b}\\
&+\frac{(a-b+k)(2n-a-b+k+2)}{8n(n+2)}\kappa.
\end{split}\nonumber
\end{equation}
Here, substituting $t=a-b\ge 0$ and $s=n-a\ge 0$, we can verify  that the coefficient of $B_{1,-b}$ is non-negative.
\end{enumerate}
On a negative quaternionic K\"ahler manifold, 
\begin{enumerate}
\item For $0\le k\le n-a$,
\begin{equation}
\begin{split}
 &dd^{\ast}+d^{\ast}d\\
=&\frac{2(a+2)(a-b+1)(n-a-k)}{(a-b+2)(k+2)(n-a)}B_{1,a+1}\\
&+\frac{2(2 n-b+5)(n-b+2)(2 k n-a k-b k+2n-2b+5k+6)
}{(a-b+2)(k+2)(2n-a-b+4)(n-b+3)}B_{1,-b}\\
&+\frac{2(k+1)(2n-a+4)(2n-a-b+3)}{(k+2)(2n-a-b+4)}B_{1,-a}\\
&+\frac{2(a+2)(a-b+1)(n-a+1)}{(a-b+2)(n-a)}B_{-1,a+1}\\
&+\frac{2(a-b+1)(2n-b+5)(n-b+2)}{(a-b+2)(2n-a-b+4)(n-b+3)}B_{-1,-b}\\
&-\frac{(a-b+k)(2n-a-b-k+2)}{8n(n+2)}\kappa.
\end{split}\nonumber
\end{equation}
\item For $n-a<k\le 2n-a-b$,
\begin{equation}
\begin{split}
&dd^{\ast}+d^{\ast}d\\
=&\frac{2(2n-b+5)(2n-a-b+3)(n-b+2)}{(a-b+2)(2n-a-b+4)(n-b+3)}B_{1,-b}\\
&+\frac{2(2n-a+4)(2n-a-b+3)(n-a+1)
}{(n-a+2)(2n-a-b+4)}B_{1,-a}\\
&+\frac{2(a+2)(a-b+1)(k+1)}{k(a-b+2)}B_{-1,a+1}\\
&+\frac{2(2n-b+5)(n-b+2)(2a+3k+a k-b k-2n)}{(a-b+2)k(2n-a-b+4)(n-b+3)}B_{-1,-b}\\
&+\frac{2(2n-a+4)(2n-a-b+3)(a+k-n)}{k(2n-a-b+4)(n-a+2)}B_{-1,-a}\\
&-\frac{(a-b+k+2)(2n-a-b-k)}{8n(n+2)}\kappa.
\end{split}\nonumber
\end{equation}
\end{enumerate}
Then we have a lower bound of the eigenvalues of $dd^{\ast}+d^{\ast}d$ on $\mathbf{S}_{k,(2_b,1_{a-b})}$ for $a>b>0$. 

From the results given in this section, we complete an eigenvalue estimate of $dd^{\ast}+dd^{\ast}$ on $S^k(\mathbf{H})\hotimes \Lambda^{a,b}_0(\mathbf{E})$. 
\begin{theorem}\label{theorem:8-2}
The eigenvalues of the Laplace operator $dd^{\ast}+d^{\ast}d$ on $S^k(\mathbf{H})\hotimes \Lambda^{a,b}_0(\mathbf{E})$ for $0\le k\le 2n-a-b$ and $0\le b\le a\le n$ have the following lower bound.
\begin{enumerate}
\item On a compact positive quaternionic K\"ahler manifold, 
\[
\begin{cases}
\displaystyle \frac{(a-b)(2n-a-b+4)}{8n(n+2)}\kappa & \textrm{for $k=0$},\\
\displaystyle \frac{(a-b+k)(2n-a-b+k+2)}{8n(n+2)}\kappa & \textrm{for $k\neq 0$}.
\end{cases}
\]
\item On a compact negative quaternionic K\"ahler manifold, 
\begin{enumerate}
\item when $a=b=0$, 
\[
-\frac{(k+2)(2n-k)}{8n(n+2)}\kappa.
\]
\item when $a=b>0$, 
\[
\begin{cases}
\displaystyle -\frac{k(2n-2a-k+4)}{8n(n+2)}\kappa & \textrm{for $0\le k\le \frac{2n-2a}{3}$},\\
\displaystyle -\frac{(k+2)(2n-2a-k)}{8n(n+2)}\kappa & \textrm{for $\frac{2n-2a}{3}< k\le 2n-2a$}.
\end{cases}
\]
\item when $a>b\ge 0$,
\[
\begin{cases}
\displaystyle -\frac{(a-b+k)(2n-a-b-k+2)}{8n(n+2)}\kappa  & \textrm{for $0\le k \le n-a$},\\
\displaystyle-\frac{(a-b+k+2)(2n-a-b-k)}{8n(n+2)}\kappa & \textrm{for $n-a< k\le 2n-a-b$}.
\end{cases}
\]
\end{enumerate}
\end{enumerate}
\end{theorem}
From this theorem, we know which irreducible bundles carry harmonic forms. The next corollary leads to the weak Lefschetz theorems for quaternionic K\"ahler manifolds in \cite{S1} and \cite{SW}.
\begin{corollary}[\cite{SW}]\label{corollary:vanish}
We consider the bundle of differential forms on a compact quaternionic K\"ahler manifold. If the scalar curvature is positive, a harmonic form is a section of $\Lambda^{a,a}_0(\mathbf{E})$ for $0\le a\le n$. If the scalar curvature is negative, a harmonic form is a section of $\Lambda^{a,a}_0(\mathbf{E})$ for $0\le a\le n$, or $S^{2n-a-b}(\mathbf{H})\hotimes \Lambda^{a,b}_0(\mathbf{E})$ for $0\le b \le a\le n$.
\end{corollary}

We shall finish by discussing relations between our eigenvalue estimates and the first eigenvalues on the quaternionic projective space $\mathbb{H}P^n$ with $\kappa=2n$. In \cite{T}, C. Tsukamoto calculated the spectra of the Laplace operator $dd^{\ast}+d^{\ast}d$ on $\mathbb{H}P^n$. On $S^k(\mathbf{H})\hotimes \Lambda^a_0(\mathbf{E})$, the first eigenvalue coincides with the lower bound in Theorem \ref{theorem:8-2}. But, so does not on $S^k(\mathbf{H})\hotimes \Lambda^{a,b}_0(\mathbf{E})$ for $a\ge b>0$. The reason is that we use only  \eqref{eqn:bw-3}--\eqref{eqn:bw-5} and \eqref{eqn:bw-6}. Since the curvature $R^{hyper}$ is zero on $\mathbb{H}P^n$, we can use all Bochner-Weitzenb\"ock formulas \eqref{eqn:bw-1}--\eqref{eqn:bw-5} to estimate eigenvalues. Then we have a better eigenvalue estimate which coincides with the first eigenvalue on $\mathbb{H}P^n$. 
\begin{example}
We consider the Laplace operator $dd^{\ast}+d^{\ast}d$ on $S^k(\mathbf{H})\hotimes \Lambda^{a,b}_0(\mathbf{E})$. We can easily show that, for $k\ge 2$, the first eigenvalue $\lambda_1$ on $\mathbb{H}P^n$ is 
\[
\frac{1}{4(n+2)}(k(k+2n+2)+a(2n-a+2)+b(2n-b+4)).
\]
On the other hand, it follows from \eqref{eqn:bw-1} and \eqref{eqn:bw-3} that 
\[
\begin{split}
dd^{\ast}+d^{\ast}d=&\sum_{N,\nu} \left(1+\frac{w_{\nu}}{2}+\frac{W_N}{2n}\right)B_{N,\nu}\\
=&\frac{1}{8(n+2)}(2k(k+2)+\pi_{(2_b,1_{a-b})}(c_2))+\sum_{N,\nu} B_{N,\nu}\\
=&\frac{1}{8(n+2)}(2k(k+2n+2)+\pi_{(2_b,1_{a-b})}(c_2))+\sum_{\nu} \frac{2(k+1)}{k+2}B_{1,\nu}\\
=&\lambda_1+\sum_{\nu} \frac{2(k+1)}{k+2}B_{1,\nu}.
\end{split}
\]
Thus we verify that the lower bound induced from Bochner-Weitzenb\"ock formulas coincides with the first eigenvalue on $\mathbb{H}P^n$. 
\end{example}

\section*{Acknowledgement}
The author was partially supported by the Grant-in-Aid for JSPS Research Fellowships for Young Scientists.


\begin{flushleft}
Yasushi Homma \\
Department of Mathematics, \\
 Faculty of Science and Technology, \\
 Tokyo University of Science, \\ 
 Noda, Chiba, 278-8510, \\
 JAPAN. \\
  \textit{E-mail address}: homma\_yasushi@ma.noda.tus.ac.jp
\end{flushleft}
\end{document}